\documentclass[12pt]{article}
\usepackage{mathrsfs}
\usepackage[pdftex]{graphicx}
\usepackage{amsmath,amssymb}
\usepackage{amsfonts}
\usepackage{txfonts}

\newtheorem{theorem}{Theorem}
\newtheorem{lemma}{Lemma}

\newcommand\figcaption{\def\@captype{figure}\caption}

\begin{document}
\title{Quasi-periodic solutions for differential equations with an elliptic-type degenerate equilibrium point under small perturbations
\thanks{This work is supported  by the
NNSF(11371132) of China, by Key Laboratory of High Performance Computing and Stochastic Information Processing. $^a$email:
lixuemei$\_1$@sina.com, $^b$ email: zaijiu@amss.ac.cn} }
\author{ Xuemei Li$^a$ and Zaijiu Shang$^b$ \\
\small $^{a}$ Key Laboratory of High Performance Computing and  Stochastic Information Processing,\\
 \small Department of Mathematics,Hunan Normal University,Changsha,Hunan 410081, P. R. China\\
\small $^{b}$  Institute of Mathematics, AMSS, Chinese Academy of Science, Beijing 100080, P. R. China.\\
}
\date{}
\maketitle \vskip 0.3in

{\bf Abstract.} This work focuses on the existence of quasi-periodic solutions for ordinary and delay differential equations (ODEs and DDEs for short) with an elliptic-type degenerate equilibrium point under quasi-periodic perturbations. We prove that under appropriate hypotheses there exist quasi-periodic solutions for perturbed ODEs and DDEs near the equilibrium point for most parameter values, then apply these results to the delayed van der Pol's oscillator with zero-Hopf singularity.

\vskip 0.2in

\renewcommand{\theequation}{\thesection.\arabic{equation}}
\section*{ 1. Introduction}
\setcounter{section}{1}\setcounter{equation}{0}

The existence problem of quasi-periodic solutions (invariant tori) is a very active research topic of KAM theory for integrable or partially integrable systems under quasi-periodic perturbations. This problem has been studied near elliptic equilibrium points by Jorba and Simo \cite{JS}, near nondegenerate invariant tori by Friedman \cite{Fr67}; Braaksma and Broer \cite{BB87}; Braaksma, Broer and Huitema \cite{BBH90}; Broer, Huitema and Sevryuk \cite{BHS96}, and Sevryuk \cite{Sev07}. Here we only mention those works closely related to the present one for non-conservative systems.

In the case where the equilibrium point or torus is degenerate (that is, the coefficient matrix of the linear part has zero as an eigenvalue, or there is a zero eigenvalue in the normal direction of the torus), the linear terms of unperturbed systems can not control the shift of equilibrium points so that this problem become complicated. Hence, the further restriction is imposed on the perturbation except the smallness and smoothness \cite{BG01, CW99, LY03, Gen07, GG05, HLY, Li16, LY12} where the averaged perturbation system in normal direction has a nondegenerate equilibrium point, or the higher-degree terms of the unperturbed system have to be taken into account \cite{You} for Hamiltonian systems and \cite{Xu} for dissipative systems.

In \cite{You, Xu}, some high-degree terms are used to control the shift of equilibrium points. You \cite{You} considered the real analytic Hamiltonian
\begin{equation}\label{hamsys}
H=\langle\omega,y\rangle+\frac{1}{2}v^2-u^{2d}+P(x,y,u,v), \qquad d\geq 2
\end{equation}
in $(x,y,u,v)$-space $\mathbb{T}^{n_0}\times \mathbb{R}^{n_0}\times \mathbb{R}^{2}$, continuously depending on the parameter $\omega$ in an open set $\mathcal{D}\subset \mathbb{R}^{n_0}$. It is proven that for beforehand fixed frequency vector $\omega_0$ satisfying the Diophantine condition
$$|\langle k,\omega_0\rangle|\geq \gamma |k|^{-\iota}, \qquad \iota >n_0+1, 0\neq k\in \mathbb{Z}^{n_0},$$
the hyperbolic-type degenerate $n_0$ dimensional torus $y=0, u=v=0$ with frequency $\omega_0$ can survive small perturbations, that is, there is an $\omega*$  such that \eqref{hamsys} at $\omega*$ possesses an n-torus with frequency $\omega_0$.

A similar problem has been investigated for degenerate lower dimensional tori for Hamiltonian systems by Cheng \cite{Ch96, Ch99}.

Concerning on dissipative systems, Xu \cite{Xu} considered the two dimensional quasi-periodic and real analytic system
\begin{equation}\label{dissys}
 \left\{\begin{array}{rl}
\dot x & =\Omega y+h_1(x,y,\omega t)+f_1(x,y,\omega t)\\
\dot y & =x^3+h_2(x,y,\omega t)+f_2(x,y,\omega t),
\end{array}
\right.
\end{equation}
where $(x,y)\in \mathbb{R}^{2}, \Omega>0$ is a constant, $h_1$ and $h_2$ are high-degree terms
$$ h_i=\sum_{n\geq 2\, {\rm or}\,l+n\geq 4}h_{iln}(\omega t)x^ly^n, \qquad i=1,2,$$
  $f_1$ and $f_2$ are lower-degree small perturbation terms
$$ f_i=\sum_{n\leq 1\, {\it and}\,l+n\leq 3}f_{iln}(\omega t)x^ly^n, \qquad i=1,2.$$
The origin (0,0) is a hyperbolic-type degenerate equilibrium point of the unperturbed system of \eqref{dissys}. In \cite{Xu}, it is proven that the system \eqref{dissys} has a quasi-periodic solution with frequency $\omega$ near the origin if the perturbation $(f_1,f_2)$ is sufficiently small and the frequency vector $\omega$ satisfies the Diophantine condition.

But the case where the degenerate equilibrium point is elliptic-type is more complicated \cite{You}: there are resonances between normal and tangent frequencies, which requires measure estimates; the linear coordinate transformation of equilibrium points at each iteration step may not smoothly depend on parameters. One aim of the present paper is to examine the existence of quasi-periodic solutions in such a case.

One the other hand, Li and Llave \cite{LL} discussed the existence of quasi-periodic solutions for delay differential equations under some assumptions, one of which is that the unperturbed linear system does not have zero as an eigenvalue. Another aim of the present paper is examine whether quasi-periodic solutions still exist for delay differential equations when the linear system has zero-eigenvalues. We will also use one high-degree term to control the shift of equilibrium point and these results will be stated in Section 2.

We conclude this section by a remark on another degeneracy in KAM theory concerning the frequency mapping of the unperturbed tori. In classical KAM theory, the frequency mapping is always required to satisfy the Kolmogorov non-degeneracy condition or some strong non-degeneracy conditions concerning the dependence on the parameters of systems. But in many concrete systems, the strong non-degeneracy conditions are not verified. This strongly motivated searching for weaker non-degeneracy conditions, which have been studied in a series of papers, for example, first by Arnol'd \cite{Ar63}, then by Bruno \cite{Bru92}, Cheng and Sun \cite{CS94}, R\"ussmann \cite{Ru01}, Sevryuk \cite{Sev07}, Han,Li and Yi \cite{HLY10} for finite dimensional Hamiltonian systems, and Bambusi, Berti and Magistrelli \cite{BBM11} for infinite dimensional case.

\renewcommand{\theequation}{\thesection.\arabic{equation}}
\section*{2. Statement of Results }
\setcounter{section}{2}\setcounter{equation}{0}

{\bf 2.1 Quasi-periodic solutions: the ODEs' case} \quad Consider the following real analytic ordinary differential equation
 \begin{equation}\label{ode1}
 \left\{\begin{array}{rl}
\dot{v}_1 & =\Omega_1(a)v_1^3+f_1(v;a)+\varepsilon g_1(\omega t,v;a,\varepsilon)\\
\dot{v}_2 & =-\Omega_2(a)v_3+d_1(a)v_1^3+f_2(v;a)+\varepsilon g_2(\omega t,v;a,\varepsilon)\\
\dot{v}_3 & =\Omega_2(a)v_2+d_2(a)v_1^3+f_3(v;a)+\varepsilon g_3(\omega t,v;a,\varepsilon),
\end{array}
\right.
 \end{equation}
where $v=(v_1,v_1,v_3)^T\in \mathbb{R}^3$ ($\cdot^T$ denotes transpose), $a\in \Pi_0\subset \mathbb{R}$ is a parameter variable, the perturbations $\varepsilon g_j (j=1,2,3)$ are quasi-periodic in $t$ with frequency $\omega=(\omega_1,\cdots,\omega_{n_0})$, $ f_j (j=1,2,3)$ are higher-degree terms
$$f_j=\sum_{l\in \Sigma^\prime}f_{j,l}(a)v^l,\qquad v^l=v_1^{l_1}v_2^{l_2}v_3^{l_3},$$
$$\Sigma^\prime= \{l=(l_1,l_2,l_3)\in \mathbb{Z}_+^3: \quad l_1\geq 4, \,{\rm or}\,l_2+l_3\geq 2,\,{\rm or}\,l_2+l_3=1\,{\rm and}\,l_1\geq 2\}.$$
The equation \eqref{ode1} is a perturbed one of the following equation with an elliptic-type degenerate equilibrium point
$$\left\{\begin{array}{rl}
\dot{v}_1 & =\Omega_1(a)v_1^3+f_1(v;a)\\
\dot{v}_2 & =-\Omega_2(a)v_3+d_1(a)v_1^3+f_2(v;a)\\
\dot{v}_3 & =\Omega_2(a)v_2+d_2(a)v_1^3+f_3(v;a).
\end{array}
\right.$$
Our aim is to prove the existence of quasi-periodic solutions of \eqref{ode1}. To this end the following assumptions are made.

\vskip 0.2in
 \qquad {\bf (H1)} All functions in \eqref{ode1} are $2\pi$-periodic in $\phi=\omega t$, which means quasi-periodic in $t$ with frequency $\omega=(\omega_1,\cdots,\omega_{n_0})$; analytic in variables $\phi\in \mathbb{T}^{n_0}, v\in \mathfrak{B}(s_0)$, $\varepsilon$ in some neighborhood of $\varepsilon=0$ (it is notable that the analyticity is not necessary but it considerably simplifies the proofs); and continuously differentiable with respect to the parameter $a\in \Pi_0$. Here $\Pi_0\subset \mathbb{R}$ is a bounded closed set of positive Lebesgue measure,  $\mathfrak{B}(s_0)=\{v\in \mathbb{R}^3: |v|\leq s_0\}$ with $s_0$ being a positive constant.

\vskip 0.2in
\qquad {\bf (H2)} The frequency vector $\omega$ satisfies the Diophantine condition
\begin{equation}\label{omegaest}
|\langle k,\omega_0\rangle|\geq \gamma_0 |k|^{-\iota},\qquad \forall k\in \mathbb{Z}^{n_0}\setminus \{0\},
\end{equation}
where $\iota \geq n_0+1$ and $0<\gamma_0\ll 1$ are constants, $|k|=|k_1|+\cdots +|k_{n_0}|$ for integer vectors, and $\langle\cdot, \cdot\rangle$ is the usual scalar product. Without loss of generality, we fix $\iota=n_0+1$ to simplify notation in proving procedure.

\vskip 0.2in
\qquad {\bf (H3)} There is a positive constant $c_0$ such that
$$\inf_{a \in \Pi_0}\left|\frac{d}{da}
\Omega_2(a)\right|\geq c_0,\qquad  \inf_{a \in \Pi_0}|\Omega_j(a)|\geq c_0,\quad j=1,2. $$

\vskip 0.2in
{\bf Remark 1.1} In this paper, the continuous differentiability of a function $f$  with respect to the parameter $a$ on a bounded closed set $\Pi$ means that the $f$ is continuous together with the first derivative in some neighborhood of $\Pi$.
\vskip 0.2in

As $g_j$ is real analytic in $v$, expand $g_j$ into power series in $v$
$$g_j=\sum_{l\in \mathbb{Z}_+^3}g_{j,l}(\omega t; a,\varepsilon)v^l, \qquad j=1,2,3. $$
Let $\phi =\omega t$,
$$g_{j,L}=\sum_{l\in \Sigma_L^0}g_{j,l}(\phi; a,\varepsilon)v^l $$
denote the lower-degree terms of $g_j$,
$$g_{j,H}=\sum_{l\in \Sigma_H^0}g_{j,l}(\phi; a,\varepsilon)v^l $$
denote the higher-degree terms of $g_j$, where
$$\Sigma_L^0=\left\{l=(l_1,l_2,l_3)\in \mathbb{Z}_+^3: \quad l_2+l_3=0\,{\rm and}\, 0\leq l_1\leq 3, \, {\rm or}\,l_2+l_3=1\,{\rm and}\, l_1=0 \right\},$$
$$\Sigma_H^0=\left\{l=(l_1,l_2,l_3)\in \mathbb{Z}_+^3:\quad l_1\geq 4,\, {\rm or}\, |l|\geq 2 \,{\rm and}\, l_2+l_3\geq 1 \right\}.$$

We want to look for quasi-periodic solutions of \eqref{ode1} with the frequency $\omega$ by a sequence of quasi-periodic transformations to eliminate the lower-degree terms $g_{j,L}$ so that \eqref{ode1} can be reduced to suitable normal forms with the coordinate origin as an equilibrium point, respectively. At each step of the transformation procedure we need to translate the coordinate origin to some equilibrium point and the change may not be smoothly but only continuously depends on the parameter even if the original equation analytically depends on the parameter. For the hyperbolic case, it is well known that we do not need to estimate the parameter measure. But for the elliptic-type degenerate case, because of resonances between the frequency $\omega$ of perturbations and the normal frequencies which requires the parameter measure estimate, we need to impose some restrictive conditions on the perturbation so that translation changes are smooth.

Define the average of the function $g_1(\phi, v; a,\varepsilon)$ in the equation \eqref{ode1} on $\phi$ with $v_2=v_3=0,\varepsilon=0$ by
 $$ \widehat{g_1}(0,v_1,0;a,0)\coloneqq \frac{1}{(2\pi)^{n_0}}\int_{\mathbb{T}^{n_0}}g_1(\phi, v_1,0,0; a,0)d\phi = \sum_{l_1\geq 0}\widehat{g_{1,l_100}}(0;a,0)v_1^{l_1}.$$

\vskip 0.2in
\qquad {\bf (H4)} Assume there is a positive constant $c_1$ such that

 \qquad {\bf Case 1} \hskip 0.8in $\inf_{a \in \Pi_0}\left|\widehat{g_{1,000}}(0;a,0)\right|\geq c_1,$\\
or

\qquad {\bf Case 2} \hskip 0.4in $\widehat{g_{1,000}}(0;a,0)= 0 \, {\rm for}\, a\in \Pi_0,\qquad {\rm and} \qquad
\inf_{a \in \Pi_0}\left|\widehat{g_{1,100}}(0;a,0)\right|\geq c_1.$

\begin{theorem}\label{theorem1}
Suppose that for the equation \eqref{ode1} Assumptions (H1)-(H4) hold. Then there is a sufficiently small $\varepsilon*>0$ such that for $0<\varepsilon\leq \varepsilon*$, there exists a Cantorian-like subset $\Pi_{\gamma_0}\subset \Pi_0$ with the Lebesgue measure
$${\rm Meas}\Pi_{\gamma_0}={\rm Meas} \Pi_0-O(\gamma_0),$$
and for any $a\in \Pi_{\gamma_0}$, the equation \eqref{ode1}
possesses a quasi-periodic solution $v=v(\omega t; a)$ which is real analytic in $\omega t$, Lipschitz in $a\in \Pi_{\gamma_0}$
 and  satisfies
 \begin{equation}\label{solest}
\sup_{\mathbb{R}\times
\Pi_{\gamma_0}}||v(\omega t; a)||=O(\varepsilon^{\frac{1}{3}}).
\end{equation}
\end{theorem}

\vskip 0.3in
{\bf 2.2 Quasi-periodic solutions: the DDEs' case} \quad Consider the perturbed delay differential equation
 \begin{equation}\label{dde0}
\dot x(t)=A(a)x(t)+B(a)x(t-1)+f(x(t),x(t-1);a)+\varepsilon g(\omega t,x(t),x(t-1);a,\varepsilon),
\end{equation}
where $x\in \mathbb{R}^q, A,B\in \mathbb{R}^{q\times q}$, the small perturbation $\varepsilon g$ is quasi-periodic in $t$ with frequency $\omega=(\omega_1,\cdots,\omega_{n_0})$, the parameter $a\in \Pi_0\subset \mathbb{R}$ and $\Pi_0$ is a bounded closed set of positive Lebesgue measure,
$$ f=O(||(x(t),x(t-1))||^3).$$
Assume the unperturbed linear equation
\begin{equation}\label{lindde}
\dot x(t)=A(a)x(t)+B(a)x(t-1)
\end{equation}
has a simple zero-eigenvalue and a pair of purely imaginary simple eigenvalues $\pm \Omega_2(a)\sqrt{-1}$, and the rest eigenvalues $\{\lambda_j(a): j=1,2,\cdots\}$ satisfy
\begin{equation}\label{rootest}
{\rm Re} \lambda_j(a)\geq \mu  \qquad {\rm for\,all} \, a\in \Pi_0,\, j=1,2,\cdots
\end{equation}
with some positive constant $\mu$.

Our another aim is to obtain the existence of quasi-periodic solutions of \eqref{dde0}. According to the center direction and hyperbolic direction--infinite dimensional part (see pages 3753-3755 in \cite{LY12} for calculation, also see the next section for a similar calculation) we can decompose the equation \eqref{dde0} into
 \begin{equation}\label{dde1}
 \left\{\begin{array}{rl}
\dot{v}_1 & =\Omega_1(a)v_1^3+f_1+\varepsilon g_1\\
\dot{v}_2 & =-\Omega_2(a)v_3+d_1(a)v_1^3+f_2+\varepsilon g_2\\
\dot{v}_3 & =\Omega_2(a)v_2+d_2(a)v_1^3+f_3+\varepsilon g_3\\
\frac{dz_t}{dt} & =U_{Q}z_t+ X_0^Q [d_3(a)v_1^3+f_4+\varepsilon g_4],
\end{array}
\right.
 \end{equation}
 where $v=(v_1,v_2,v_3)^T$, $U_Q$ is the restriction of $U$ to the hyperbolic invariant subspace $Q$ of \eqref{lindde}, $U$ is the extension to $C^1([-1,0],\mathbb{R}^q)$ of the infinitesimal generator $U_0$ of the solution operator semigroup of the linear equation \eqref{lindde}, $z_t\in Q, X_0^Q$ is the projection of $X_0$ on the hyperbolic invariant subspace $Q$ of \eqref{lindde}, $ X_0(\theta)=0$ if $ -1\leq \theta<0$;  $=E_q$ if
$\theta=0$, and $E_q$ is the $q\times q$ identity matrix, and
  for $j=1,\cdots,4$,
 $$f_j=f_j(v,z_t(0),z_t(-1);a)=\sum_{(l,m,n)\in \Sigma^{\prime\prime}}f_{j,lmn}(a)v^lz_t^m(0)z_t^n(-1),$$
 $$ g_j=g_j(\omega t, v,z_t(0),z_t(-1);a,\varepsilon)=\sum_{l,m,n}g_{j,lmn}(\omega t; a,\varepsilon)v^lz_t^m(0)z_t^n(-1) $$
 with $(l,m,n)\in \mathbb{Z}_+^3\times \mathbb{Z}_+^q\times \mathbb{Z}_+^q $, and
$$\Sigma^{\prime\prime}=\left\{(l,m,n)\in \mathbb{Z}_+^3\times \mathbb{Z}_+^q\times \mathbb{Z}_+^q: \, l_1\geq 4, \,{\rm or}\,|l|+|m|+|n|\geq 3\,{\rm and}\, l_2+l_3+|m|+|n|\geq 1\right\}.$$
 The spectral set of $U_Q$ is just $\{\lambda_j(a): j=1,2,\cdots\}$ satisfying \eqref{rootest}.

If removing the hyperbolic direction from \eqref{dde1}, then the equation \eqref{dde1} is reduced to a similar form of \eqref{ode1}, and we also have an analogous result on the existence of quasi-periodic solutions of \eqref{dde1}.

\vskip 0.2in
\qquad {\bf (H1)$^\prime$}Assume all functions in \eqref{dde1} are $2\pi$-periodic in $\phi=\omega t$, analytic in variables $\phi\in \mathbb{T}^{n_0}, v$ and $z_t$ on $\mathfrak{B}(s_0)$, $\varepsilon$ in some neighborhood of $\varepsilon=0$, and continuously differentiable with respect to the parameter $a$ on $\Pi_0$, where
$$\mathfrak{B}(s_0)\coloneqq \{(v,z)\in \mathbb{R}^3\times Q:\quad |v|\leq s_0, ||z||\leq s_0\},$$
$z\in Q\subset C([-1,0],\mathbb{R}^q)$ and $||z||=\sup_{-1\leq \theta\leq 0}|z(\theta)|$, $s_0$ is a positive constant.
\vskip 0.2in

Define the average of the function $g_1(\phi, v, z_t(0), z_t(-1); a,\varepsilon)$ in the equation \eqref{dde1}  on $\phi$ with $v_2=v_3=0, z_t=0, \varepsilon=0$ by
$$
 \widehat{g_1}(0,v_1,0;a,0) \coloneqq \frac{1}{(2\pi)^{n_0}}\int_{\mathbb{T}^{n_0}}g_1(\phi, v_1,0,0,0,0; a,0)d\phi
   = \sum_{l_1 \geq 0}\widehat{g_{1,(l_100)00}}(0;a,0)v_1^{l_1}.$$

\qquad {\bf (H4)$^\prime$} Assume there is a positive constant $c_1$ such that

 \qquad {\bf Case 1} \hskip 0.8in $\inf_{a \in \Pi_0}\left|\widehat{g_{1,(000)00}}(0;a,0)\right|\geq c_1,$\\
or

\qquad {\bf Case 2} \hskip 0.4in $\widehat{g_{1,(000)00}}(0;a,0)= 0 \, {\rm for}\, a\in \Pi_0,\quad {\rm and} \quad
\inf_{a \in \Pi_0}\left|\widehat{g_{1,(100)00}}(0;a,0)\right|\geq c_1.$

\begin{theorem}\label{theorem2}
Suppose that for the equation \eqref{dde1} Assumptions (H1)$^\prime$, (H2)-(H3) and (H4)$^\prime$ hold, and all eigenvalues of the operator $U_Q$ satisfy the condition \eqref{rootest}. Then there is a sufficiently small $\varepsilon*>0$ such that for $0<\varepsilon\leq \varepsilon*$, the conclusions in Theorem 1 appear to hold true also for the equation \eqref{dde1}. Correspondingly, the equation \eqref{dde0} possesses a quasi-periodic solution for any $ a\in\Pi_{\gamma_0}$.
\end{theorem}

In Section 3, we will apply Theorem 2 to delayed van der Pol's oscillator with zero-Hopf singularity, and obtain the existence of quasi-periodic solutions of the system under quasi-periodic perturbations. The proofs of our results are based on a rapidly convergent iteration process and reducing the equations to some normal forms with zero as an equilibrium point. The sections 4 and 5 are devoted to the proof of Theorem 1. In Section 4, we consider one step of the iteration and obtain several estimates, that is, we give an iteration lemma and its proof, which shows that at each iteration step the transformed equation possesses much smaller lower-degree terms. The convergence of the iteration sequence of quasi-periodic transformations and the measure estimate of parameter sets are analysed in Section 5. As the hyperbolic direction in (\ref{dde1}), involving the infinite dimensional part, does not result in essential difficulties and its treatment is similar to that in \cite{LY12}, and the argument for the center direction is the same as in (\ref{ode1}), we will only describe the proof sketch of Theorem 2 in Section 6.

\renewcommand{\theequation}{\thesection.\arabic{equation}}
\section*{3. Quasi-periodic response in delayed van der Pol's oscillator}
\setcounter{section}{3}\setcounter{equation}{0}

Consider the delayed van der Pol's oscillator
\begin{equation}\label{dvdpe}
 \ddot x(t)+a(x^2(t)-1)\dot x(t)+x(t) = f(x(t-\tau))+\varepsilon g(\omega^\prime t,x(t),x(t-\tau),\varepsilon),
 \end{equation}
where $a$ is a constant, the delay $\tau\geq 0$, $f$ and $g$ are assumed to be analytic in all variables, $\varepsilon g$ is a quasi-periodic perturbation in time $t$ with the frequency $\omega^\prime =(\omega_1^\prime,\cdots,\omega_{n_0}^\prime)$.

Braaksma and Broer \cite{BB87}, Broer, Huitema and Sevryuk \cite{BHS96} investigated the existence of quasi-periodic solutions of (\ref{dvdpe}) for the case of the delay $\tau=0$.

For the case without the perturbation term, the bifurcation and stability of (\ref{dvdpe}) have been extensively studied, see \cite{BDL14, JW08, WJ10, XC03, ZG13} and the references therein.

Assume
$$ f(0)=f^{\prime\prime}(0)=0, \qquad f^\prime(0)=b, \qquad f^{\prime\prime\prime}(0)=6b_1$$
just as in \cite{WJ10, XC03, ZG13}. The characteristic equation of the linearization of unperturbed equation of (\ref{dvdpe}) at the equilibrium point $x=0$ is
\begin{equation}\label{chareq}
\lambda^2-a\lambda+1-b e^{-\tau \lambda}=0.
\end{equation}
The distribution of roots of (\ref{chareq}) is analyzed by Jiang and Wei \cite{JW08}. Especially, in the case where $b=1, \tau=\tau_0(a)$ and $0<a<\sqrt 2$, the characteristic equation (\ref{chareq}) has a simple zero root and a pair of purely imaginary simple roots $\pm \omega_0\sqrt{-1}$, and the rest roots of (\ref{chareq}) have negative real parts (see Lemma 2.7 in \cite{JW08}), where
\begin{equation}\label{OT}
\omega_0=\sqrt{2-a^2},\qquad   \tau_0= \left\{\begin{array}{ll}
 \frac{1}{\omega_0}\arcsin a\omega_0,\quad 1<a<\sqrt 2\\
  \frac{1}{\omega_0}(\pi -\arcsin a\omega_0),\quad 0<a\leq 1.
 \end{array} \right.
 \end{equation}

We will discuss the existence of quasi-periodic solutions of Equation (\ref{dvdpe}) in this case. We will regard $a$ as a parameter variable, and take a closed subset in its allowed domain, such as $\Pi_0=[\frac{1}{4},\frac{5}{4}]$ for simplicity.

We rewrite Equation (\ref{dvdpe}) with $b=1,\tau=\tau_0$ (meanwhile, rescaling the  time $t$ by $\tau_0 t$) as
\begin{equation}\label{rdvdpe}
\dot y(t)=Ay(t) +By(t-1)+F_1(y(t),y(t-1))+F_2(\tau_0\omega^{\prime}t, y(t),y(t-1),\varepsilon),
\end{equation}
where
\begin{eqnarray*}
 y(t) & = & {\rm col}(x(t), \dot x(t)),\\
 F_1 & = & {\rm col}(0, \tau_0(-ax^2(t)\dot x(t)+b_1x^3(t-1))),\\
 F_2 & = & {\rm col}(0, \tau_0[f(x(t-1))-x(t-1)-b_1x^3(t-1)+\varepsilon g(\tau_0\omega^{\prime}t, x(t),x(t-1),\varepsilon)]),
 \end{eqnarray*}
 $$A = \left(\begin{array}{cc} 0 & \tau_0\\ -\tau_0 & a\tau_0\end{array}\right), \qquad B=\left(\begin{array}{cc} 0 & 0\\  \tau_0 & 0\end{array}\right).$$

Because in the considered case, the linear system of (\ref{rdvdpe}) has a zero eigenvalue, we will regard $F_1$ ana $F_2$ as a principal term and a perturbation, respectively. We need to decompose the equation (\ref{rdvdpe}) into the center and hyperbolic directions. Following an analogous procedure in Section 1 in \cite{LY12}, the phase space of (\ref{rdvdpe}) is $\mathcal{C}\coloneqq C([-1,0],\mathbb{R}^2)$ endowed with the supremum norm. Let $U_0$ be the infinitesimal generator of the solution operator semigroup for the linearization system of (\ref{rdvdpe}) with $\varepsilon=0$, given by
$$
 U_0z=\frac{dz}{d\theta} \qquad {\rm for}\,
z\in \mathscr{D}(U_0)\coloneqq \left\{z\in \mathcal{C}:\quad  \frac{dz}{d\theta}\in
\mathcal{C}, \frac{dz}{d\theta}|_{\theta=0}=Az(0)+Bz(-1)\right\},$$
and $U$ be its extension to
$\mathcal{C}^1\coloneqq C^1([-1,0],\mathbb{R}^2)$ defined by
$$U z  \coloneqq
\frac{dz}{d\theta}+X_0\left[Az(0)+Bz(-1)-\frac{dz}{d\theta}|_{\theta=0}\right], \qquad U:\mathcal{C}^1\rightarrow \mathcal{BC}\coloneqq \mathcal{C}\oplus \langle X_0\rangle,$$
where $ X_0(\theta)=0$ if $ -1\leq \theta<0$;  $=E_2$ if
$\theta=0$, and $E_2$ is the $2\times 2$ identity matrix, $\langle X_0\rangle=\{X_0(\theta)c: c\in \mathbb{R}^2\}$, $\mathcal{BC}$ is the enlarged phase space of all bounded functions from $[-1,0]$ to $\mathbb{R}^2$, continuous on $[-1,0)$, with a possible jump discontinuity at $\theta=0$, which can be identified with $\mathcal{C}\times \mathbb{R}^2$. The operators $U_0$ and $U$ have only point spectrum consisting of all roots multiplied by $\tau_0$, of the characteristic equation (\ref{chareq}) with $\tau=\tau_0$.

Let  $\mathcal {C}^{\ast} \coloneqq C([0,1],\mathbb{R}^{2
^{\ast}})$, where $\mathbb{R}^{2 ^{\ast}}$ is the $2$-dimensional
row-vector space. Define a
bilinear form $\langle \cdot,\cdot\rangle$ on
$\mathcal{C}^{\ast}\times \mathcal{C}$ by
$$\langle \psi,z\rangle= \psi(0)z(0)+\int_{0}^{1}\psi(s)Bz(s-1)ds,\qquad \psi\in
\mathcal{C}^{\ast}, z\in \mathcal{C}, $$
and the formal adjoint
operator $U_0^{\ast}$ of $U_0$ by
$$\langle \psi,U_0z\rangle=\langle U_0^{\ast}\psi,z\rangle \qquad
{\it for} \quad z\in \mathscr{D}(U_0)\quad {\it and }\quad \psi \in
\mathscr{D}(U_0^{\ast}).$$
Then
$$ U_0^{\ast}\psi=-\frac{d\psi}{ds} \qquad {\rm for}\, \psi \in
\mathscr{D}(U_0^{\ast})= \left\{\psi\in \mathcal{C}^{\ast}:\quad
\frac{d\psi}{ds}\in \mathcal{C}^{\ast},
\frac{d\psi}{ds}|_{s=0}=-\psi(0)A -\psi(1)B \right\}.$$
Set
$$\Lambda=\{0,\pm\tau_0\omega_{0}\sqrt{-1}\},$$
which is the set of eigenvalues of $U$ with zero real parts. Thus, there exists a positive constant  $\mu$, independent of $a$, such that
\begin{equation}\label{erpc}
{\rm Re} \lambda(a)\leq -\mu
\end{equation}
for $a\in \Pi_0$, where the $\lambda(a)$  denotes an arbitrary eigenvalue of $U$ with a nonzero real part.

The phase space $\mathcal{C}$ is decomposed by $\Lambda$ as
$\mathcal{C}=P_{\Lambda}\oplus Q_{\Lambda}$, and $\Phi_{\Lambda}=(z_1(\theta),z_2(\theta),z_3(\theta))(\theta \in [-1,0])$ is a real basis for the generalized eigenspace $P_{\Lambda}$, where
$$z_1(\theta)=\left(\begin{array}{c}1 \\ 0\end{array}\right),\qquad
z_2(\theta)=\left(\begin{array}{c}\sin(\tau_0\omega_0\theta) \\ \omega_0\cos(\tau_0\omega_0\theta) \end{array}\right),\qquad
z_3(\theta)=\left(\begin{array}{c}\cos(\tau_0\omega_0\theta) \\ -\omega_0\sin(\tau_0\omega_0\theta) \end{array}\right).$$
The daul basis $\Psi_{\Lambda}$ of $\Phi_{\Lambda}$ (i.e., the real basis for the generalized eigenspace of $U_0^{\ast}$ with respect to $\Lambda$), satisfying
 $\langle
\Psi_{\Lambda},\Phi_{\Lambda} \rangle=E_3$ (the $3\times3$ identity matrix), is $\Psi_{\Lambda}(s)={\rm col}(\psi_1(s),\psi_2(s),\psi_3(s)) (s\in [0,1])$,
\begin{eqnarray*}
 \psi_1(s) & = & (\frac{a}{a-\tau_0},-\frac{1}{a-\tau_0}),\\
 \psi_2(s) & = & \frac{1}{2(b_2^2+b_3^2)}(\omega_0(\tau_0-a)\cos(\tau_0\omega_0s)+(4-a^2-a\tau_0)\sin(\tau_0\omega_0s), \\
  & &\hskip 0.9in -2b_2\sin(\tau_0\omega_0s)-2b_3\cos(\tau_0\omega_0s)),\\
 \psi_3(s) & = & \frac{1}{2(b_2^2+b_3^2)}((4-a^2-a\tau_0)\cos(\tau_0\omega_0s)-\omega_0(\tau_0-a)\sin(\tau_0\omega_0s),\\
 & &  \hskip 0.9in 2b_3\sin(\tau_0\omega_0s)-2b_2\cos(\tau_0\omega_0s)),
\end{eqnarray*}
where
$$b_2=\frac{1}{2}(a-\tau_0\cos(\tau_0\omega_0)),\qquad b_3=-\omega_0+\frac{\tau_0}{2}\sin(\tau_0\omega_0).$$
From the characteristic equation (\ref{chareq}), it follows
$$b_2=\frac{1}{2}(a+\tau_0-\tau_0a^2),\qquad b_3=\frac{\omega_0}{2}(a\tau_0-2).$$
Therefore,
$$P_{\Lambda}=\{z\in \mathcal{C}:\quad z=\Phi_{\Lambda}b_0 \quad
{\rm for \; some}\quad b_0\in \mathbb{R}^3\},$$
$$Q_{\Lambda}=\{z\in \mathcal{C}:\quad \langle \Psi_{\Lambda},z \rangle=0 \}, \qquad U\Phi_{\Lambda}=\Phi_{\Lambda}U_{\Lambda},$$
where $U_{\Lambda}$ is a $3\times 3$ matrix
$$U_{\Lambda}=\left(\begin{array}{ccc} 0 & 0 & 0\\ 0 & 0 & -\tau_0 \omega_0\\ 0 & \tau_0 \omega_0 & 0 \end{array}\right).$$
If $y(t)$ is a solution of (\ref{rdvdpe}), then $y_t$ (defined by $y_t(\theta)=y(t+\theta), -1\leq \theta\leq 0$) may be decomposed as
$$y_t=\Phi_{\Lambda}\langle \Psi_{\Lambda},y_t \rangle + y_t^Q \coloneqq \Phi_{\Lambda}v(t)+z_t,$$
where $v(t)\in \mathbb{R}^3, z_t\in Q_{\Lambda}$, but it does not necessarily satisfy $z_t(\theta)=z(t+\theta)$ for $-1\leq \theta \leq 0$. Especially, we have
$$y(t)=\Phi_{\Lambda}(0)v(t)+z_t(0),\qquad y(t-1)=\Phi_{\Lambda}(-1)v(t)+z_t(-1).$$
Let $U_Q$ denote $U$ restricted to $Q_{\Lambda}$. The equation (\ref{rdvdpe}) can be written as

\begin{equation}\label{dvdpe1}
 \left\{\begin{array}{rl}
\dot{v}(t) & =U_{\Lambda}v(t)+ \Psi_{\Lambda}(0)F(\tau_0\omega^{\prime}t, y(t),y(t-1),\varepsilon)\\
\frac{dz_t}{dt} & =U_{Q}z_t+ X_0^QF(\tau_0\omega^{\prime}t, y(t),y(t-1),\varepsilon),
\end{array}
\right.
 \end{equation}
 where $F=F_1+F_2, X_0^Q(\theta)=X_0(\theta)-\Phi_{\Lambda}(\theta)\Psi_{\Lambda}(0), -1\leq \theta \leq 0$.

 Let $v=(v_1,v_2,v_3)^T, z_t=(z_{1t},z_{2t})^T$.  Noting
 \begin{equation}\label{xexpr}
 x(t)=v_1(t)+v_3(t)+z_{1t}(0),  \qquad x(t-1)=v_1(t)-v_2(t)\sin(\tau_0\omega_0)+v_3(t)\cos(\tau_0\omega_0)+z_{1t}(-1),
 \end{equation}
and employing the Taylor series expansion of $f$, the equation (\ref{dvdpe1}) reads
 \begin{equation}\label{dvdpe2}
 \left\{\begin{array}{rl}
\dot{v}_1 & =-\frac{\tau_0b_1}{a-\tau_0}v_1^3-\frac{\tau_0}{a-\tau_0}\mathfrak{G}\\
\dot{v}_2 & =-\tau_0\omega_0v_3-\frac{\tau_0b_3}{b_2^2+b_3^2}[b_1v_1^3+\mathfrak{G}]\\
\dot{v}_3 & =\tau_0\omega_0v_2-\frac{\tau_0b_2}{b_2^2+b_3^2}[b_1v_1^3+\mathfrak{G}]\\
\frac{dz_t}{dt} & =U_{Q}z_t+ X_0^Q \cdot {\rm col} (0, b_1v_1^3+\mathfrak{G}),
\end{array}
\right.
 \end{equation}
 where
 $$\mathfrak{G}=\varepsilon g(\tau_0\omega^{\prime}t, x(t),x(t-1),\varepsilon)+\sum\limits_{(l,n)\in\Sigma^{\prime\prime}}f_{ln}v^lz_{1t}^{n_1}(0)z_{2t}^{n_2}(0)z_{1t}^{n_3}(-1),$$
 $$v^l=v_1^{l_1}v_2^{l_2}v_3^{l_3}, \qquad |l|=l_1+l_2+l_3,\qquad |n|=n_1+n_2+n_3,$$
 $$\Sigma^{\prime\prime}=\left\{(l,n)\in \mathbb{Z}_+^3\times \mathbb{Z}_+^3: \, l_1\geq 4, \,{\rm or}\,|l|+|n|\geq 3\,{\rm and}\, l_2+l_3+|n|\geq 1\right\}.$$

 We obtain quasi-periodic solutions of (\ref{dvdpe}) by employing Theorem 2 in the case where the characteristic equation \eqref{chareq} has a simple zero root and a pair of purely imaginary roots (that is, $b=1$ and $\tau=\tau_0$ defined by (\ref{OT})). Thus we also want to impose a restrictive condition on the lower-degree terms of $g$. Expand $g$ as the power series in $(x(t),x(t-\tau))$
\begin{equation}\label{gexp}
g(\omega^\prime t, x(t), x(t-\tau); \varepsilon)=\sum_{m+n\geq 0}g_{mn}(\omega^\prime t;\varepsilon)x^m(t)x^n(t-\tau).
\end{equation}
\vskip 0.2in
\qquad {\bf (H4)$^{\prime\prime}$} Assume that there is at least one whose average does not vanish, of the following two functions $g_{00}(\phi;0)$ and $g_{10}(\phi;0)+g_{01}(\phi;0)$ in $\phi=\omega^\prime t \in \mathbb{T}^{n_0}$ from coefficients of the lower-degree terms in (\ref{gexp}).

\begin{theorem}\label{theorem3}
Suppose that for the equation (\ref{dvdpe}) with $b=1,\tau=\tau_0$ and $\omega^\prime$ satisfying (\ref{omegaest}), Assumption (H4)$^{\prime\prime}$ holds and $b_1\neq 0$. Then for most parameter values $a$ in the set $\Pi_0=[\frac{1}{4},\frac{5}{4}]$ and for sufficiently small $\varepsilon$, the equation (\ref{dvdpe}) possesses a quasi-periodic solution $x(t)$ with frequency $\omega^\prime$ and $x=O(\varepsilon^{\frac{1}{3}})$.
\end{theorem}

{\bf proof}\quad We prove Theorem 3 by employing Theorem 2, which means we only need to verify Assumptions (H1)$^\prime$,(H3) and (H4)$^\prime$ for the equation \eqref{dvdpe2}. In fact, using the expression \eqref{xexpr} we rewrite the equation \eqref{dvdpe2} in the form of \eqref{dde1}, and Assumption (H1)$^\prime$ holds by the assumption on the equation \eqref{dvdpe}. Moreover,
$$\Omega_1(a)=-\frac{\tau_0b_1}{a-\tau_0}, \qquad \Omega_2(a)=\tau_0\omega_0$$
and
\begin{equation}\label{4g1expr}
g_1(\phi,v,z_t(0),z_t(-1);a,\varepsilon)=-\frac{\tau_0}{a-\tau_0}g(\phi,v_1+v_3+z_{1t}(0),v_1-v_2\sin(\tau_0\omega_0)+v_3\cos(\tau_0\omega_0)+z_{1t}(-1); \varepsilon),
\end{equation}
where $\phi=\tau_0\omega^\prime t$. The expression \eqref{OT} and $b_1\neq 0$ imply (H3) with $a\in \Pi_0=[\frac{1}{4},\frac{5}{4}]$.

Noting that from the expression \eqref{4g1expr} it follows
$$\widehat{g_1}(0,v_1,0;a,0)=-\frac{\tau_0}{a-\tau_0}\sum_{m+n\geq 0}\widehat{g_{mn}}(0;0)v_1^{m+n},$$
the condition (H4)$^{\prime\prime}$ implies (H4)$^{\prime}$, which completes the proof of Theorem 3.

\renewcommand{\theequation}{\thesection.\arabic{equation}}
\section*{4. The iteration step}
\setcounter{section}{4}\setcounter{equation}{0}

Introducing complex conjugate coordinates, which are more convenient to use in the following, into the second and third equations of (\ref{ode1}), that is
\begin{equation}\label{complexT}
w_1=v_1,\qquad w_2=v_2+\sqrt{-1}v_3,\qquad \bar w_2=v_2-\sqrt{-1}v_3,
\end{equation}
the equation (\ref{ode1}) reads
 \begin{equation}\label{ode2}
 \left\{\begin{array}{rl}
\dot{w}_1 & =\Omega_1(a)w_1^3+F_1(w;a)+\varepsilon G_1(\phi,w;a,\varepsilon)\\
\dot{w}_2 & =\sqrt{-1}\Omega_2(a)w_2+d_4(a)w_1^3+F_2(w;a)+\varepsilon G_2(\phi,w;a,\varepsilon)\\
\dot{\bar w}_2 & =-\sqrt{-1}\Omega_2(a)\bar w_2+\overline{d_4(a)}w_1^3+F_3(w;a)+\varepsilon G_3(\phi,w;a,\varepsilon),
\end{array}
\right.
 \end{equation}
where $\phi=\omega t$, $w=(w_1,w_2,\bar w_2)^T$, $d_4(a)=d_1(a)+\sqrt{-1}d_2(a)$, $F_j$ and $G_j$, expressed in $f_i$ and $g_i (i=1,2,3)$, are the same forms as $f_j$ and $g_j$ ($j=1,2,3$) respectively, and satisfy the following reality condition (see $\S 15$ in \cite{SM71})
\vskip 0.2in
\qquad {\bf (H5) (Reality Condition)}
$$ \overline{H_1}(\phi,w_1,\bar w_2,w_2;a,\varepsilon)=H_1(\phi,w_1,w_2,\bar w_2;a,\varepsilon),$$
$$ \overline{H_2}(\phi,w_1,\bar w_2,w_2;a,\varepsilon)=H_3(\phi,w_1,w_2,\bar w_2;a,\varepsilon),$$
$$\overline{H_3}(\phi,w_1,\bar w_2,w_2;a,\varepsilon)=H_2(\phi,w_1,w_2,\bar w_2;a,\varepsilon),$$
where $H_1, H_2$ and $H_3$ are the right-hand-side expressions (RHSE, for short) of $\dot w_1, \dot w_2$ and $\dot{\bar w}_2$  in (\ref{ode2}), respectively. Here the complex conjugation $\bar F$ of a function $F$ means complex conjugation of the coefficients in the power series of $F$.
Noting $\widehat{G_1}(0,w_1,0,0;a,\varepsilon)=\widehat{g_1}(0,w_1,0,0;a,\varepsilon)$, the assumption (H4) implies
\begin{equation}\label{case1}
\inf_{a \in \Pi_0}\left|\widehat{G_{1,000}}(0;a,0)\right|\geq c_1,
\end{equation}
or
\begin{equation}\label{case2}
\widehat{G_{1,000}}(0;a,0)=0 \quad {\rm for}\, a\in \Pi_0,\qquad {\rm and} \qquad
\inf_{a \in \Pi_0}\left|\widehat{G_{1,100}}(0;a,0)\right|\geq c_1.
\end{equation}

In order to regard the terms $d_4(a)w_1^3$ and $\overline{d_4(a)}w_1^3$ in the second and third equations of (\ref{ode2}) as small perturbations, rescaling
\begin{equation}\label{rescaleT}
w_1\rightarrow \varepsilon^{\frac{1}{4}}w_1,\qquad w_2\rightarrow \varepsilon^{\frac{1}{2}}w_2,\qquad \varepsilon_0=\varepsilon^{\frac{1}{4}},
\end{equation}
 the equation (\ref{ode2}) is changed into the following form
\begin{equation}\label{ode3}
 \left\{\begin{array}{rl}
\dot{w}_1 & =\varepsilon_0^2\left[\Omega_1(a)w_1^3+\varepsilon_0 F_1^0(\phi,w;a,\varepsilon_0)+\sum_{l\in \mathfrak{\Sigma}_L^0}\varepsilon_0^{2|l|+1-l_1} G_{1,l}(\phi; a,\varepsilon_0^4)w^l\right] \\
\dot{w}_2 & =\sqrt{-1}\Omega_2(a)w_2+\varepsilon_0 d_4(a)w_1^3+\varepsilon_0^2F_2^0(\phi,w;a,\varepsilon_0)+ \sum_{l\in \mathfrak{\Sigma}_L^0}\varepsilon_0^{2|l|+2-l_1} G_{2,l}(\phi; a,\varepsilon_0^4)w^l\\
\dot{\bar w}_2 & =-\sqrt{-1}\Omega_2(a)\bar w_2+\varepsilon_0 \overline{d_4(a)}w_1^3+\varepsilon_0^2F_3^0(\phi,w;a,\varepsilon_0)+ \sum_{l\in \mathfrak{\Sigma}_L^0}\varepsilon_0^{2|l|+2-l_1} G_{3,l}(\phi; a,\varepsilon_0^4)w^l  ,
\end{array}
\right.
 \end{equation}
where, $F_j^0,j=1,2,3,$ are the higher-degree terms.

As an analytic function in real variables can always extend to a real analytic function on a complex domain of its variables, by the assumption (H1) there are positive constants $r_0,s_0$ and $M_0$ such that all of terms in (\ref{ode3}) are analytic in $(\phi,w)\in \mathscr{D}(r_0,s_0)\coloneqq \mathscr{U}(r_0)\times \mathscr{W}(s_0)$, continuously differentiable with respect to the parameter $a\in \Pi_0$, satisfying Reality Condition and
$$|\|F_j^0|\|_{r_0,s_0,\Pi_0}\coloneqq \max_{s=0,1}\sup_{\mathscr{D}(r_0,s_0)\times
 \Pi_0}|\partial_{a}^s F_j^0|\leq M_0,\qquad j=1,2,3,$$
 where $\partial_a$ denotes the partial derivative with respect to $a$, $\mathscr{U}(r_0)$ and $\mathscr{W}(s_0)$ are complex neighborhoods of the torus $\mathbb{T}^{n_0}$ and the origin respectively,
 $$\mathscr{U}(r_0)=\{\phi \in  \mathbb{C}^{n_0}/2\pi\mathbb{Z}^{n_0}:\, |{\rm Im}\phi|\coloneqq \sup_{1\leq j\leq n_0}|{\rm Im}\phi_j| \leq r_0\},\qquad
 \mathscr{W}(s_0)=\{w \in  \mathbb{C}^3:\, |w| \leq s_0\}.$$

 To obtain quasi-periodic solutions of (\ref{ode3}), the main idea is to make the lower-degree terms smaller and smaller by a sequence of coordinate transformations. So we assume that at the general $\nu$-th step it is true, then we want to look for a coordinate transformation which is defined in a slightly smaller domain in $\phi=\omega t$, such that in the new coordinate (\ref{ode3}) has much smaller lower-degree terms than the ones in the $\nu$-th step.

 Before stating the iteration lemma, we list sequences for the pertinent parameters and notations. For $\nu\geq 0$, set\\
(i) $\epsilon_{\nu+1}=\epsilon_{\nu}^{\frac{9}{8}} (\nu\geq 1),\quad \epsilon_1=\varepsilon_0^{\iota},\quad \epsilon_0=\varepsilon_0^4; \qquad \gamma_{\nu}=\gamma_0/(\nu+1)^2, $\\
 the $\iota>0$ is a constant;\\
(ii) $K_\nu=-\frac{1}{r_0}(\nu+1)^22^{\nu+2}\ln \epsilon_1,\quad K_0=-\frac{8}{r_0}\ln \varepsilon_0;$\\
(iii) $r_\nu=(1-\sigma_\nu)r_0,\quad  s_\nu=(1-\sigma_\nu)s_0$ with
$\sigma_0=0, \sigma_\nu=(1^{-2}+\cdots+\nu^{-2})/(2\sum_{j=1}^{\infty}j^{-2}),$ \\
 $\rho_\nu=r_\nu-r_{\nu+1},\quad
\delta_\nu=\frac{1}{3}(s_\nu-s_{\nu+1}),\quad \mathscr{W}_\nu\coloneqq \mathscr{W}(s_\nu),\quad \mathscr{D}_\nu\coloneqq \mathscr{D}(r_\nu,s_\nu)=\mathscr{U}(r_\nu)\times \mathscr{W}(s_\nu)$;\\
(iv) For two quantities $u$ and $v$, if there is an absolute
 constant $C$ such that $|u|\leq C|v|(|u|\geq C|v|)$ with
 some metric $|\cdot|$, we write $u\lessdot v (u\gtrdot v)$;\\
(v) $f\in \mathfrak{A}_{r,s}^{\Pi}$ means that $f(\phi,x;a)$ is analytic in
 $(\phi,x)\in \mathscr{D}(r,s)$, and
 continuously differentiable in $a\in \Pi$; moreover, for a given $\epsilon$, the $f$ satisfies
 $$|\|f|\|_{r,s,\Pi}\coloneqq \max_{j=0,1}||\partial_{a}^j f||_{r,s,\Pi}\lessdot \epsilon,\qquad || f||_{r,s,\Pi}\coloneqq \sup_{\mathscr{D}(r,s)\times
 \Pi}| f|,$$
 denoted by $f=O_{r,s,\Pi}(\epsilon)$;\\
(vi)  for $f\in \mathfrak{A}_r^{\Pi}$
 with Fourier series expansion
 $$f(\phi;a)=\sum_{k\in \mathbb{Z}^{n_0}}\hat
 f(k;a)e^{\sqrt{-1}\langle k,\phi\rangle},\qquad \phi\in \mathscr{U}(r),a\in \Pi,$$
 define the truncation operator $\Gamma_K$ as follows
 $$\Gamma_K f=\sum_{|k|\leq K}\hat
 f(k;a)e^{\sqrt{-1}\langle k,\phi\rangle},$$
 where $|k|=|k_1|+\cdots+|k_{n_0}|$. We have
 \begin{equation}\label{fcoeffe}
 ||\partial_a^j\hat f(k;a)||_{\Pi}\leq
 |\|f|\|_{r,\Pi}e^{-|k|r},\qquad j=0,1
 \end{equation}
 and
 \begin{equation}\label{trune}
||\partial_a^j({\rm Id}-\Gamma_K)f||_{r-\rho,\Pi}\leq
 C|\|f|\|_{r,\Pi}\rho^{-n_0}e^{-\rho K},\qquad j=0,1, 0<\rho <r.
 \end{equation}

 \begin{lemma}\label{iterativelemma}
 Assume $\epsilon_1\leq \gamma_0^6$ and that at the $\nu$-th ($\nu\geq 1$) step we have obtained the following quasi-periodic system
\begin{equation}\label{odev}
 \left\{\begin{array}{rl}
\dot{w}_1 & =\varepsilon_0^2\left[N_1^\nu(w_1;a)+ F_1^\nu(\phi,w;a)+G_1^\nu(\phi,w; a)\right] \\
\dot{w}_2 & =\Omega_2^\nu(a)w_2+ F_2^\nu(\phi,w;a)+G_2^\nu(\phi,w; a)\\
\dot{\bar w}_2 & =\Omega_3^\nu(a)\bar w_2+ F_3^\nu(\phi,w;a)+G_3^\nu(\phi,w; a)
\end{array}
\right.
 \end{equation}
with the normal form
$$N_1^\nu(w_1;a)=\Omega_1^\nu(a)w_1^3+e_2^\nu(a)w_1^2+e_1^\nu(a)w_1,$$
satisfying\\
{\bf $(\nu. 1)$ (Frequency condition)} \qquad $\overline{\Omega_2^\nu}=\Omega_3^\nu,$ and for $j=1,2$
$$\inf_{a\in \Pi_{\nu-1}}|\partial_a({\rm Im}\Omega_2^\nu(a))|\geq c_0(1-\sigma_\nu),\quad \inf_{a\in\Pi_{\nu-1}}|\Omega_j^\nu(a)|\geq c_0(1-\sigma_\nu),$$
$$\inf_{a\in \Pi_{\nu-1}}|e_1^\nu(a))|\geq c_2(1-\sigma_\nu)\varepsilon_0^\kappa, \quad |\|e_1^\nu|\|_{\Pi_{\nu-1}}\lessdot \varepsilon_0^\kappa,
\quad |\|e_2^\nu|\|_{\Pi_{\nu-1}}\lessdot \varepsilon_0^\beta, $$
$$|\|\Omega_j^\nu(a)-\Omega_j^{\nu-1}(a)|\|_{\Pi_{\nu-1}}\lessdot\epsilon_{\nu-1},\quad |\|e_j^\nu(a)-e_j^{\nu-1}(a)|\|_{\Pi_{\nu-1}}\lessdot\epsilon_{\nu-1}^{\frac{4}{3}}(\nu \geq 2), $$
where $c_2, \kappa$ and $\beta$ are positive constants, and $\beta <\frac{4}{3}\iota$, $\kappa<\frac{4}{3}\iota+\beta$, $e_j^0=0 $;\\
{\bf $(\nu. 2)$ (Smallness condition)}  for $j=1,2,3,$ $F_j^\nu, G_j^\nu\in \mathfrak{A}_{r_\nu,s_\nu}^{\Pi_\nu}$,
$$F_j^\nu=\sum_{l\in \mathfrak{\Sigma}_H^0}F^\nu_{j,l}(\phi; a)w^l$$
are the higher-degree terms, and
$$G_j^\nu=\sum_{l\in \mathfrak{\Sigma}_L^0}G^\nu_{j,l}(\phi; a)w^l$$
are the lower-degree small perturbation terms with
\begin{equation}\label{Fj}
|\|F_j^\nu-F_j^{\nu-1}|\|_{r_\nu,s_\nu,\Pi_\nu}\leq \epsilon_\nu^{\frac{1}{3}},\qquad    |\|F_j^\nu|\|_{r_\nu,s_\nu,\Pi_\nu}\leq M_0(1+\sigma_\nu),
 \end{equation}
\begin{equation}\label{Gl1}
 G^\nu_{j,l_100}=O_{r_\nu,\Pi_\nu}(\epsilon_\nu^{4-l_1}),\qquad 0\leq l_1\leq 3,
 \end{equation}
\begin{equation}\label{Gl2l3}
 G^\nu_{j,0l_2l_3}=O_{r_\nu,\Pi_\nu}(\epsilon_\nu),\qquad l_2+l_3=1;
 \end{equation}
{\bf $(\nu. 3)$ (Reality condition)}  denoting $w=(w_1,w_2,\bar w_2)^T, \underline{w}=(w_1,\bar w_2,w_2)^T$, the RHSEs of $\dot w_1, \dot w_2$ and $\dot {\bar w}_2$ in (\ref{odev}) by $H_1^\nu, H_2^\nu$ and $H_3^\nu$, respectively, for $\phi \in \mathbb{T}^{n_0}$
$$ \overline{H_1^\nu}(\phi,\underline{w};a)=H_1^\nu(\phi,w;a),\quad \overline{H_2^\nu}(\phi,\underline{w};a)=H_3^\nu(\phi,w;a),\quad \overline{H_3^\nu}(\phi,\underline{w};a)=H_2^\nu(\phi,w;a),$$
where $\phi=\omega t$ with $\omega$ satisfying Assumption (H2).

Then there is a closed subset $\Pi_{\nu+1}\subset \Pi_\nu$ of the measure estimate
\begin{equation}\label{measest}
 {\rm Meas} \Pi_{\nu+1}\geq {\rm Meas}\Pi_\nu(1-c_3\gamma_\nu),
 \end{equation}
with a constant $c_3>0$, and a quasi-periodic coordinate transformation
$\mathscr{T}^{\nu}:\mathscr{W}_{\nu+1}\rightarrow \mathscr{W}(s_\nu-\delta_\nu)\subset \mathscr{W}_{\nu}$ in the form
$$\mathscr{T}^{\nu}:\qquad w=w_+ +W_0^\nu +W^\nu(\phi, w_+; a),$$
where $\phi\in \mathscr{U}(r_{\nu+1}), a\in \Pi_{\nu+1}$, $w_+=(w_{1+},w_{2+},\bar w_{2+})^T$ is the new coordinate, $W_0^\nu=(w_{10}^\nu(a),0,0)^T$ with
$ w_{10}^\nu(a)\in \mathbb{R}$, $W^\nu=(W^\nu_1,W^\nu_2,W^\nu_3)^T, W^\nu_j \in \mathfrak{A}_{r_{\nu+1},s_{\nu+1}}^{\Pi_{\nu+1}}$ is of the same form as $G_j^\nu(j=1,2,3)$, consisting of lower-degree terms, and satisfies Reality condition
\begin{equation}\label{rc}
\overline{W_1^\nu}(\phi,\underline{w_+};a)=W_1^\nu(\phi,w_+;a),\quad \overline{W_2^\nu}(\phi,\underline{w_+};a)=W_3^\nu(\phi,w_+;a),\quad \overline{W_3^\nu}(\phi,\underline{w_+};a)=W_2^\nu(\phi,w_+;a)
 \end{equation}
for $\phi \in \mathbb{T}^{n_0}$, and the following estimates
\begin{eqnarray}
w_{10}^\nu & = & O_{\Pi_{\nu+1}}(\epsilon_{\nu}^{\frac{4}{3}}), \label{w10est}\\
W_{1,l_100}^\nu & = & O_{r_{\nu+1},\Pi_{\nu
+1}}(\varepsilon_0^2\epsilon_{\nu}^{4-l_1}\gamma_{\nu}^{-2}\rho_{\nu}^{-(3n_0+2)}),\qquad 0\leq l_1\leq 3,\label{W1l1}\\
W_{1,0l_2l_3}^\nu & = & O_{r_{\nu+1},\Pi_{\nu
+1}}(\varepsilon_0^2\epsilon_{\nu}\gamma_{\nu}^{-2}\rho_{\nu}^{-(3n_0+2)}),\qquad l_2+l_3=1,\label{W1l2l3}\\
W_{j,l_100}^\nu & = & O_{r_{\nu+1},\Pi_{\nu
+1}}(\epsilon_{\nu}^{4-l_1}\gamma_{\nu}^{-2}\rho_{\nu}^{-(3n_0+2)}),\qquad 0\leq l_1\leq 3; j=2,3,\label{Wjl1}\\
W_{j,0l_2l_3}^\nu & = & O_{r_{\nu+1},\Pi_{\nu
+1}}(\epsilon_{\nu}\gamma_{\nu}^{-2}\rho_{\nu}^{-(3n_0+2)}),\qquad l_2+l_3=1, j=2,3,\label{Wjl2l3}
\end{eqnarray}
such that the equation (\ref{odev}) is transformed into the same form as (\ref{odev}) satisfying Conditions $(\nu.1)-(\nu.3)$ by replacing $\nu$ with $\nu+1$ and $(w_1,w_2,\bar w_2)$ with $(w_{1+},w_{2+},\bar w_{2+})$, respectively. Here we have suppressed the dependence of all functions on
$\varepsilon_0$.
\end{lemma}

{\bf proof} In the proof of the lemma we drop the index $\nu$ and write `+' for $`\nu +1'$ to
 simplify  notation. Thus, $F_j=F_j^\nu, F_j^+=F_j^{\nu+1}$, and so on. Also, we drop the parameter $a$ from functions whenever there is no confusion.
 \vskip 0.1in
 {\bf 1) Quasi-periodic transformation}

 Set
$$N(w)={\rm col}(N_1(w_1),\Omega_2 w_2,\Omega_3 \bar w_2),\qquad F(\phi,w)={\rm col}(F_1(\phi,w),F_2(\phi,w),F_3(\phi,w)),$$
$$G(\phi,w)={\rm col}(G_1(\phi,w),G_2(\phi,w),G_3(\phi,w)),\qquad J_{\varepsilon_0}={\rm diag}(\varepsilon_0^2,1,1).$$
Then we can rewrite the equation (\ref{odev}) as
\begin{equation}\label{rodev}
J_{\varepsilon_0}^{-1}\dot w=N(w)+ G(\phi,w)+F(\phi,w).
 \end{equation}

We will obtain the transformation $\mathscr{T}$ by two steps: the first step is to find a quasi-periodic transformation $\mathscr{T}_1$ to reduce the lower-degree terms $G$ into a simplest form with constant coefficients, and the second step is to translate the equilibrium point to the origin. Take the  $\mathscr{T}_1$ as follows
$$ \mathscr{T}_1:\qquad w=u+U(\phi,u)$$
with $u=(u_1,u_2,\bar u_2)^T$ and
$$U_j(\phi,u)=\sum_{l\in \mathfrak{\Sigma}_L^0}U_{j,l}(\phi)u_1^{l_1}u_2^{l_2}\bar u_2^{l_3},\qquad j=1,2,3.$$
Substituting the transformation $\mathscr{T}_1$ into (\ref{rodev}), we have
\begin{eqnarray}
& &J_{\varepsilon_0}^{-1}(E_3 +DU)J_{\varepsilon_0}(J_{\varepsilon_0}^{-1}\dot u-N(u))= - J_{\varepsilon_0}^{-1}\omega\partial_{\phi}U-J_{\varepsilon_0}^{-1}DU\cdot J_{\varepsilon_0}N(u)+DN(u)\cdot U+\Gamma_K G(\phi,u) \nonumber\\
 & &\quad +({\rm Id}-\Gamma_K)G(\phi,u)+R_1(\phi,u)+\int_0^1DG(\phi,u+\xi U)\cdot U d\xi+F\circ \mathscr{T}_1,\label{nodev}
 \end{eqnarray}
where
\begin{equation}\label{R1}
R_1(\phi,u)={\rm col}((3\Omega_1 u_1+e_2)U_1^2+\Omega_1U_1^3, 0,0),
 \end{equation}
$\partial_\phi U$ denotes the partial derivative of $U$ with respect to $\phi$, $DU$ and $DG$ are the Jacobian matrices of $U$ and $G$ with respect to $u$ and $w$, respectively. The latter line consists of higher-degree terms and smaller lower-degree terms. Hence, the transformation $\mathscr{T}_1$ for which we want to look satisfies
\begin{equation}\label{Ueq}
J_{\varepsilon_0}DN(u)\cdot U+J_{\varepsilon_0} \Gamma_K G(\phi,u)-\omega\partial_{\phi}U-DU\cdot J_{\varepsilon_0}N(u)=J_{\varepsilon_0}\widehat{N}(u)+J_{\varepsilon_0}R_2,
 \end{equation}
where $\widehat{N}$ is the drift from $G$ and can not be removed by the transformation $\mathscr{T}_1$, $R_2$ consists of higher-degree terms.

Inserting the polynomial expressions of all functions in \eqref{Ueq}, we obtain the following equations
\begin{eqnarray}
 & & \omega \partial_\phi U_{1,l_100}+(l_1-1)\varepsilon_0^2e_1U_{1,l_100} + (l_1-3)\varepsilon_0^2e_2U_{1,(l_1-1)00}\nonumber\\
& & \quad +(l_1-5)\varepsilon_0^2\Omega_1U_{1,(l_1-2)00}  = \varepsilon_0^2(\Gamma_K G_{1,l_100}(\phi)-\widehat{G_{1,l_100}}(0)),\quad 0\leq l_1\leq 3,\label{U1l1eq}\\
& & \omega \partial_\phi U_{1,0l_2l_3}-\varepsilon_0^2e_1U_{1,0l_2l_3} +  l_2\Omega_2U_{1,0l_2l_3} +l_3\Omega_3 U_{1,0l_2l_3}
 = \varepsilon_0^2\Gamma_K G_{1,0l_2l_3}(\phi), l_2+l_3=1,\label{U1l2l3eq}\\
& &  \omega \partial_\phi U_{j,l_100}-\Omega_j U_{j,l_100}+l_1\varepsilon_0^2e_1U_{j,l_100} + (l_1-1)\varepsilon_0^2e_2U_{j,(l_1-1)00} \nonumber\\
 & & \quad +(l_1-2)\varepsilon_0^2\Omega_1U_{j,(l_1-2)00} = \Gamma_K G_{j,l_100}(\phi),\hskip 0.5in 0\leq l_1\leq 3; j=2,3,\label{U23l1eq}
\end{eqnarray}
\begin{equation}\label{U2l2l3eq}
\omega \partial_\phi U_{2,010}=\Gamma_K G_{2,010}(\phi)-\widehat{G_{2,010}}(0),\quad \omega \partial_\phi U_{2,001}+(\Omega_3-\Omega_2)U_{2,001}= \Gamma_K G_{2,001}(\phi),
 \end{equation}
\begin{equation}\label{U3l2l3eq}
\omega \partial_\phi U_{3,001}=\Gamma_K G_{3,001}(\phi)-\widehat{G_{3,001}}(0),\quad \omega \partial_\phi U_{3,010}-(\Omega_3-\Omega_2)U_{3,010}= \Gamma_K G_{3,010}(\phi),
 \end{equation}
and the higher-degree terms $R_2={\rm col}(R_{21},\varepsilon_0^2R_{22},\varepsilon_0^2R_{23})$ with
\begin{eqnarray}
 R_{21} & = & (\Omega_2 U_{1,200}-e_2 U_{1,300})u_1^4+\sum_{l_2+l_3=1}(2e_2+3\Omega_1 u_1)U_{1,0l_2l_3}u_1u_2^{l_2}\bar u_2^{l_3},\label{R21}\\
 R_{22} & = & -[2\Omega_1 U_{2,200}+3(\Omega_1 u_1+e_2)U_{2,300}]u_1^4,\label{R22}\\
  R_{23} & = & -[2\Omega_1 U_{3,200}+3(\Omega_1 u_1+e_2)U_{3,300}]u_1^4,\label{R23}
 \end{eqnarray}
here, we let $U_{j,l_100}=0$ if $l_1<0$ for $j\in \{1,2,3\}$. We take
\begin{equation}\label{Nhat}
\widehat{N}(u)={\rm col}\left(\sum_{l_1=0}^3 \widehat{G_{1,l_100}}(0)u_1^{l_1}, \widehat{G_{2,010}}(0)u_2,\widehat{G_{3,001}}(0)\bar u_2\right).
 \end{equation}
The reality condition $(\nu.3)$ implies $ \widehat{G_{1,l_100}}(0),l_1=0,1,2,3,$ are real and
$$ \overline{\widehat{G_{2,010}}(0)}=\widehat{G_{3,001}}(0).$$

Set
$$\Omega_1^+=\Omega_1+\widehat{G_{1,300}}(0),\quad \Omega_2^+=\Omega_2+\widehat{G_{2,010}}(0),\quad \Omega_3^+=\Omega_3+\widehat{G_{3,001}}(0).$$
Then Conditions $(\nu.1)$ and $(\nu.2)$ imply $\Omega_j^+(j=1,2,3) $ satisfy requirements for them in Condition $(\nu.1)$ replacing $\nu$ by $\nu+1$.

We solve the equations \eqref{U1l1eq} and \eqref{U23l1eq} according to the order $l_1=0,1,2$ and 3, and each equation of \eqref{U1l1eq}-\eqref{U3l2l3eq} by a standard procedure in KAM theory. Expanding all functions of $\phi$ into Fourier series and comparing coefficients, we can find there are the following small divisors
$$\sqrt{-1}\langle k,\omega\rangle+(j-1)\varepsilon_0^2e_1,\quad \sqrt{-1}\langle k,\omega\rangle+j\varepsilon_0^2e_1-\Omega_{2,3},\qquad, 0\leq j \leq 3,$$
$$\sqrt{-1}\langle k,\omega\rangle+ \Omega_{2,3}-\varepsilon_0^2e_1,\quad \sqrt{-1}\langle k,\omega\rangle\pm (\Omega_2-\Omega_3).$$
These equations are solvable if we remove some parameter values from $\Pi$ so that small divisors satisfy the Melnikov conditions, and the Fourier series of unknown functions are convergent. Noting that the frequency $\omega$ satisfies \eqref{omegaest} and the reality condition implies $e_1$ is real, thus, we take the subset $\Pi_+$ of $\Pi$ as follows
$$
 \Pi_+=\left\{a\in \Pi: \, | \sqrt{-1}\langle k,\omega\rangle+m_2\Omega_2+m_3\Omega_3|\geq \frac{\gamma}{|k|^{n_0+1}}, 0<|k|\leq K, |m_2+m_3|\leq 1,1\leq |m_2|+|m_3|\leq 2\right\},
 $$
where $k\in \mathbb{Z}^{n_0}, m_i\in \mathbb{Z},i=2,3,$ and obtain the measure estimate \eqref{measest} and solution estimates of \eqref{U1l1eq} by the proof of Lemma 2 in \cite{Mos67} (using \eqref{fcoeffe})
\begin{eqnarray*}
||U_{1,000}||_{r_+,\Pi_+} & \lessdot & \varepsilon_0^2\gamma^{-1}\rho^{-(2n_0+1)}|\|G_{1,000}|\|_{r,\Pi},\\
||\partial_a U_{1,000}||_{r_+,\Pi_+} & \lessdot & \varepsilon_0^2\gamma^{-2}\rho^{-(3n_0+2)}|\|G_{1,000}|\|_{r,\Pi},
\end{eqnarray*}
that is,
$$|\| U_{1,000}|\|_{r_+,\Pi_+}  \lessdot \varepsilon_0^2 \gamma^{-2}\rho^{-(3n_0+2)}|\|G_{1,000}|\|_{r,\Pi}.$$
Similarly,
\begin{eqnarray*}
|\|U_{1,100}|\|_{r_+,\Pi_+} & \lessdot & \varepsilon_0^2 \gamma^{-1}\rho^{-(2n_0+1)}(|\|G_{1,100}|\|_{r,\Pi}+\varepsilon_0^2\gamma^{-2}\rho^{-2(n_0+1)}|\|G_{1,000}|\|_{r,\Pi}),\\
|\| U_{1,j00}|\|_{r_+,\Pi_+} & \lessdot & \varepsilon_0^2\gamma^{-2}\rho^{-(3n_0+2)}\left[|\|G_{1,j00}|\|_{r,\Pi} +\varepsilon_0^2\gamma^{-1}\rho^{-(n_0+1)}|\|G_{1,100}|\|_{r,\Pi}\right.\\
& & \quad \left.+\varepsilon_0^4\gamma^{-2}\rho^{-2(n_0+1)}|\|G_{1,000}|\|_{r,\Pi}\right],\quad j=2,3.
\end{eqnarray*}
Noting that if $\epsilon_1$ is sufficiently small, then $\gamma^{-2}\rho^{-2(n_0+1)}\epsilon\lessdot 1$, the condition \eqref{Gl1} and the above inequalities
imply
\begin{equation}\label{U1l1est}
|\| U_{1,l_100}|\|_{r_+,\Pi_+}  \lessdot \varepsilon_0^2\gamma^{-2}\rho^{-(3n_0+2)}\epsilon^{4-l_1},\qquad 0\leq l_1\leq 3.
 \end{equation}
More detail is seen in \cite{LL, Yuan02} for the measure estimate and continuous differentiability in $a\in \Pi_+$, and in \cite{LY12} for estimates of solutions to \eqref{U1l1eq}-\eqref{U3l2l3eq} and omitted. We obtain the following estimates of solutions to \eqref{U1l2l3eq}-\eqref{U3l2l3eq}
\begin{eqnarray}
|\|U_{1,0l_2l_3}|\|_{r_+,\Pi_+} & \lessdot & \varepsilon_0^2\gamma^{-2}\rho^{-(3n_0+2)}\epsilon,\qquad l_2+l_3=1, \label{U1l2l3est}\\
|\|U_{j,l_100}|\|_{r_+,\Pi_+} & \lessdot & \gamma^{-2}\rho^{-(3n_0+2)}\epsilon^{4-l_1},\qquad 0\leq l_1\leq 3;j=2,3, \label{Ujl1est}\\
|\|U_{j,0l_2l_3}|\|_{r_+,\Pi_+} & \lessdot & \gamma^{-2}\rho^{-(3n_0+2)}\epsilon,\quad l_2+l_3=1;j=2,3. \label{Ujl2l3est}
\end{eqnarray}
Obviously, $\mathscr{T}_1:\,\mathscr{W}(s-2\delta)\rightarrow \mathscr{W}(s-\delta)$ if $\epsilon_1$ is sufficiently small. And $U_j\in \mathfrak{A}_{r_+,s-2\delta}^{\Pi_+}$, and satisfies Reality condition like $G_j, j=1,2,3$.
 \vskip 0.2in
 {\bf 2) Translation}

Let
$$\widetilde{N}(u)=N(u)+\widehat{N}(u).$$
We want to change $\widetilde{N}$ to the same form as $N$ by a translation, which means we only need to translate the first coordinate $u_1$. Consider the equilibrium point equation
\begin{equation}\label{epeq}
\Omega_1^+u_1^3+(e_2+\widehat{G_{1,200}}(0))u_1^2+(e_1+\widehat{G_{1,100}}(0))u_1+\widehat{G_{1,000}}(0)=0.
 \end{equation}
From $||\widehat{G_{1,l_100}}(0)||_{\Pi}\lessdot \epsilon ^{4-l_1}(0\leq l_1\leq 3)$ and Condition $(\nu.1)$ it follows
$$\inf_{a\in \Pi}|e_1+\widehat{G_{1,100}}(0)|\geq \frac{c_2}{2}\varepsilon_0^\kappa.$$
By Lemma~\ref{root} (see Appendix), the equation \eqref{epeq} has a read root $w_{10}(a)$ satisfying
\begin{equation}\label{u10est}
|\|w_{10}|\|_{\Pi}\lessdot \epsilon^{\frac{4}{3}}.
 \end{equation}
Take the translation
$$\mathscr{T}_2:\qquad u_1=w_{10}+w_{1+},\quad u_2=w_{2+},\quad \bar u_2=\bar w_{2+}$$
and
$$N^+(w_+)=\widetilde{N}\circ \mathscr{T}_2(w_+)={\rm col}(N_1^+(w_{1+}),\Omega_2^+w_{2+},\Omega_3^+\bar w_{2+}).$$
Then, for sufficiently small $\epsilon_1$ it is clear that $\mathscr{T}_2:\,\mathscr{W}(s-3\delta)\rightarrow \mathscr{W}(s-2\delta)$ and
$$N_1^+(w_{1+})=\Omega_1^+w_{1+}^3+e_2^+w_{1+}^2+e_1^+w_{1+}$$
with
\begin{eqnarray*}
e_2^+-e_2 & = & 3\Omega_1^+ w_{10}+\widehat{G_{1,200}}(0),\\
e_1^+-e_1 & = & 3\Omega_1^+ w_{10}^2+2(e_2+\widehat{G_{1,200}}(0))w_{10}+\widehat{G_{1,100}}(0).
\end{eqnarray*}
Thus, by \eqref{u10est} and $||\widehat{G_{1,l_100}}(0)||_{\Pi}\lessdot \epsilon ^{4-l_1}(0\leq l_1\leq 3)$, we obtain the estimates in $(\nu.1)$ replacing $\nu$ by $\nu+1$.

Take $\mathscr{T}=\mathscr{T}_1\circ \mathscr{T}_2$. Then $\mathscr{T}:\,\mathscr{W}(s_+)\rightarrow \mathscr{W}(s-\delta)$ and the transformation $\mathscr{T}$ is defined by
\begin{eqnarray*}
w_1 & = & w_{1+}+w_{10}+U_1(\phi,w_{1+}+w_{10},w_{2+},\bar w_{2+})\coloneqq w_{1+}+w_{10} +W_1(\phi,w_{1+},w_{2+},\bar w_{2+}),\\
w_j & = & w_{j+}+U_j(\phi,w_{1+}+w_{10},w_{2+},\bar w_{2+})\coloneqq w_{j+}+W_j(\phi,w_{1+},w_{2+},\bar w_{2+}),\quad j=1,2.
\end{eqnarray*}
As $w_{10}$ is real and $U_j$ satisfies Reality condition, it implies that $W_j$ also satisfies Reality condition \eqref{rc} and $W_j\in \mathfrak{A}_{r_+,s_+}^{\Pi_+}(j=1,2,3)$. By \eqref{U1l1est}-\eqref{Ujl2l3est} and \eqref{u10est} we have the estimates \eqref{w10est}-\eqref{Wjl2l3}, and
\begin{equation}\label{Jacobiest}
|\|D\mathscr{T}|\|_{r_+,s_+,\Pi_+} \leq 1+\epsilon^{\frac{1}{3}}.
\end{equation}
 \vskip 0.2in
 {\bf 3) Estimates of remainder terms}

Inserting the translation $\mathscr{T}_2$ into the equation \eqref{nodev} and using \eqref{Ueq}, we have
\begin{equation}\label{odev+1}
J_{\varepsilon_0}^{-1}\dot w_+=N^+(w_+)+J_{\varepsilon_0}^{-1}(E_3 + DU\circ \mathscr{T}_2)^{-1}J_{\varepsilon_0}[P^1+P^2+P^3],
 \end{equation}
where
$$P^1=- \left(J_{\varepsilon_0}^{-1}DU\cdot J_{\varepsilon_0}\widehat{N}\right)\circ \mathscr{T}_2+(R_1+R_2)\circ \mathscr{T}_2+ \int_0^1(DG(\phi,{\rm Id}+\xi U)\cdot U)\circ \mathscr{T}_2 d\xi,$$
$$ P^2=({\rm Id}-\Gamma_K) G\circ \mathscr{T}_2, \qquad  P^3=F\circ \mathscr{T}.$$

We have already proved $N^+(w_+)$ satisfies the condition $(\nu.1)$ with $\nu+1$. Now, our aim is to rewrite the equation \eqref{odev+1} in the form of \eqref{odev} and prove the conditions $(\nu.2)$ and $(\nu.3)$ hold for \eqref{odev+1} replacing $\nu$ by $\nu+1$. By Lemma~\ref{reality}, Condition $(\nu.3)$ and \eqref{rc}, it is easy to see that the equation \eqref{odev+1} satisfies the reality condition $(\nu.3)$ with $\nu+1$.  In the following, we verify the condition $(\nu.2)$ with $\nu+1$ for the equation \eqref{odev+1}.

Noting that the $P^1$ is a polynomial in $w_+$, the part $P_H^1$ of its higher-degree terms satisfies
\begin{equation}\label{P1H}
|\|P_H^1|\|_{r_+,s_+,\Pi_+}\lessdot \gamma^{-2}\rho^{-(3n_0+2)}\epsilon.
 \end{equation}
Let
$$P^i=(P_1^i,P_2^i,P_3^i)^T,\qquad i=1,2,3$$
and introduce the notation
$$|\|P|\|_{s_1,s_2,s_3,\Pi_+}\coloneqq \max_{j=0,1}\sup\{|\partial_a^jP(\phi,w_{1+},w_{2+},\bar w_{2+};a)|:\, |{\rm Im} \phi|\leq s_1,|w_{1+}|\leq s_2, |w_{2+}|\leq s_3, |\bar w_{2+}|\leq s_3,a\in \Pi_+\}.$$
The estimates \eqref{U1l1est}-\eqref{Ujl2l3est} and \eqref{u10est}, the expressions \eqref{R1},\eqref{R21}-\eqref{R23} and \eqref{Nhat}, and the condition $(\nu.2)$ imply
$$|\|P_j^1|\|_{r_+,\epsilon,\epsilon^2,\Pi_+}\lessdot \gamma^{-2}\rho^{-(3n_0+2)}\epsilon^5,\qquad j=1,2,3,$$
which, by the the Cauchy inequality, implies
\begin{equation}\label{P1jl1}
|\|P_{j,l_100}^1|\|_{r_+,\Pi_+}\lessdot \gamma^{-2}\rho^{-(3n_0+2)}\epsilon^{5-l_1},\qquad 0\leq l_1\leq 3,
 \end{equation}
\begin{equation}\label{P1jl2l3}
|\|P_{j,0l_2l_3}^1|\|_{r_+,\Pi_+}\lessdot \gamma^{-2}\rho^{-(3n_0+2)}\epsilon^2,\qquad l_2+l_3=1; j=1,2,3.
 \end{equation}
The $P^2$ consists only of lower-degree terms, and the definition of the truncation operator $\Gamma_K$ (using \eqref{trune}) and \eqref{u10est} lead to
 $$|\|P_{j,l_1l_2l_3}^2|\|_{r_+,\Pi_+}\lessdot \rho^{-n_0}e^{-\rho K} |\|G_{j,l_1l_2l_3}|\|_{r,\Pi},$$
which, by the definition of $K_{\nu}$ and \eqref{Gl1}-\eqref{Gl2l3}, implies
\begin{equation}\label{P2j}
|\|P_{j,l_100}^2|\|_{r_+,\Pi_+}\lessdot \rho^{-n_0}\epsilon^{5-l_1},\qquad |\|P_{j,0l_2l_3}^2|\|_{r_+,\Pi_+}\lessdot \rho^{-n_0}\epsilon^2,\quad 0\leq l_1\leq 3, l_2+l_3=1; j=1,2,3.
 \end{equation}

To estimate $P^3$, we divide $F$ into two parts $F=F^\prime+F^{\prime\prime}$,
$$F_j^\prime(\phi,w)=F_{j,400}(\phi)w_1^4+F_{j,110}(\phi)w_1w_2+F_{j,101}(\phi)w_1\bar w_2,\qquad F_j^{\prime\prime}=F_j-F_j^\prime, \, j=1,2,3.$$
By the Cauchy inequality and \eqref{Fj}, we have
$$|\|F_{j,400}|\|_{r,\Pi}\leq \frac{M_0}{s^4}(1+\sigma),\qquad |\|F_{j,1l_2l_3}|\|_{r,\Pi}\leq \frac{M_0}{s^2}(1+\sigma),\quad l_2+l_3=1,$$
\begin{equation}\label{Fjp}
|\|F_j^{\prime\prime}|\|_{r,s,\Pi}\lessdot 4M_0(1+\sigma),
 \end{equation}
and
\begin{equation}\label{Fje}
|\|F_j\circ \mathscr{T}-F_j|\|_{r_+,s_+,\Pi_+}\leq \delta^{-1}|\|F_j|\|_{r,s,\Pi} |\|\mathscr{T}-{\rm Id}|\|_{r_+,s_+,\Pi_+}
 \lessdot \delta^{-1} \gamma^{-2}\rho^{-(3n_0+2)}\epsilon |\|F_j|\|_{r,s,\Pi}.
 \end{equation}
Let
$$P^{3\prime}(\phi,w_+)=F^\prime(\phi,\mathscr{T}(w_+)),\qquad P^{3\prime\prime}(\phi,w_+)=F^{\prime\prime}(\phi,\mathscr{T}(w_+)).$$
Then $P^3=P^{3\prime}+P^{3\prime\prime}$, and the estimates \eqref{U1l1est}-\eqref{Ujl2l3est} and \eqref{u10est} imply
\begin{eqnarray}
|\|P_{j,000}^{3\prime}|\|_{r_+,\Pi_+} & \lessdot &  \gamma^{-2}\rho^{-(3n_0+2)}\epsilon^{\frac{16}{3}},\label{P3j0}\\
|\|P_{j,100}^{3\prime}|\|_{r_+,\Pi_+} & \lessdot &  \gamma^{-2}\rho^{-(3n_0+2)}\epsilon^4,\label{P3j1}\\
|\|P_{j,200}^{3\prime}|\|_{r_+,\Pi_+} & \lessdot &  \left(\epsilon^{\frac{8}{3}}+\gamma^{-2}\rho^{-(3n_0+2)}\epsilon^3\right),\label{P3j2}\\
|\|P_{j,300}^{3\prime}|\|_{r_+,\Pi_+} & \lessdot &  \left(\epsilon^{\frac{4}{3}}+\gamma^{-2}\rho^{-(3n_0+2)}\epsilon^2\right),\label{P3j3}\\
|\|P_{j,0l_2l_3}^{3\prime}|\|_{r_+,\Pi_+} & \lessdot &  \epsilon^{\frac{4}{3}},\qquad l_2+l_3=1 \label{P3j4}
\end{eqnarray}
for $j=1,2,3$. Using the Cauchy inequality again and \eqref{Fjp} we obtain
$$|\|P_j^{3\prime\prime}|\|_{r_+,\epsilon,\epsilon^2,\Pi_+}\lessdot \epsilon^5,$$
hence,
\begin{eqnarray}
|\|P_{j,l_100}^{3\prime\prime}|\|_{r_+,\Pi_+} & \lessdot &  \epsilon^{5-l_1},\qquad 0\leq l_1\leq 3,\label{P3j5}\\
|\|P_{j,0l_2l_3}^{3\prime\prime}|\|_{r_+,\Pi_+} & \lessdot &  \epsilon^2,\qquad l_2+l_3=1; \quad j=1,2,3. \label{P3j6}
\end{eqnarray}
Here, we use the fact that $s\geq s_0/2$ and $s_0$ is a constant.

Let $P_{jL}^3$ and $P_{jH}^3$ be the lower-degree terms and higher-degree terms of $P_j^3$, respectively. Then the estimates \eqref{P3j0}-\eqref{P3j6} lead to
$$|\|P_{jL}^3|\|_{r_+,s_+,\Pi_+}\lessdot \left(\epsilon^{\frac{4}{3}}+\gamma^{-2}\rho^{-(3n_0+2)}\epsilon^2\right)\leq \epsilon_+$$
and
\begin{equation}\label{P3H}
|\| P_{jH}^3 |\|_{r_+,s_+,\Pi_+}\leq |\| P_j^3-F_j |\|_{r_+,s_+,\Pi_+}+|\|F_j |\|_{r,s,\Pi}+|\|P_{jL}^3|\|_{r_+,s_+,\Pi_+}
\leq M_0(1+\sigma)(1+\epsilon_+^{\frac{1}{3}})
 \end{equation}
by using \eqref{Fje} for $j=1,2,3$.

Set
$$P=P^1+P^2+P^3\equiv (P_1,P_2,P_3)^T.$$
Thus, by \eqref{P1jl1}-\eqref{P2j} and \eqref{P3j0}-\eqref{P3j6}, we have
\begin{equation}\label{Pjest}
|\|P_{j,l_100}|\|_{r_+,\Pi_+} < \epsilon_+^{4-l_1},\quad  |\|P_{j,0l_2l_3}|\|_{r_+,\Pi_+} < \epsilon_+,\qquad   0\leq l_1\leq 3, l_2+l_3=1; j=1,2,3
 \end{equation}
for sufficiently small $\epsilon_1$.

Set
$$\mathscr{P}=J_{\varepsilon_0}^{-1}(E_3 + DU\circ \mathscr{T}_2)^{-1}J_{\varepsilon_0}P.$$
Then
$$G^+=\mathscr{P}_L,\qquad F^+=\mathscr{P}_H,$$
where $\mathscr{P}_L$ and $  \mathscr{P}_H$ are the lower-degree terms and higher-degree terms of $\mathscr{P}$, respectively. As the estimates \eqref{U1l1est}-\eqref{Ujl2l3est} and \eqref{u10est} imply
$$J_{\varepsilon_0}^{-1}(E_3 + DU\circ \mathscr{T}_2)^{-1}J_{\varepsilon_0}=E_3+O(\gamma^{-2}\rho^{-(3n_0+2)}\epsilon),$$
hence
$$|\|J_{\varepsilon_0}^{-1}(E_3 + DU\circ \mathscr{T}_2)^{-1}J_{\varepsilon_0}|\|_{r_+,s_+,\Pi_+}\leq 1+O(\epsilon_+^{\frac{1}{3}}).$$
From \eqref{P1H}, \eqref{Fje}, \eqref{P3H} and \eqref{Pjest}, it easily follows that $G^+$ and $F^+$ satisfy the condition $(\nu.2)$ with $\nu+1$. The proof of the lemma is complete.
\vskip 0.2in

\renewcommand{\theequation}{\thesection.\arabic{equation}}
\section*{5. Proof of Theorem 1 }
\setcounter{section}{5}\setcounter{equation}{0}

We first reduce the equation \eqref{ode3} into the form of \eqref{odev}, then use Lemma~\ref{iterativelemma} to prove Theorem 1. The reducing procedure is similar to the one in the proof of Lemma~\ref{iterativelemma}, only difference is the translation. Just as doing in the proof of Lemma~\ref{iterativelemma}, introducing the change of coordinates
$$ \mathscr{T}_1^0:\quad w=u+U(\phi,u),\quad U_j(\phi,u)=\sum_{l\in \mathfrak{\Sigma}_L^0}U_{j,l}(\phi)u_1^{l_1}u_2^{l_2}\bar u_2^{l_3},\, j=1,2,3,$$
substituting the change into \eqref{ode3}, and noting that in the present case, the normal form is $N(w)={\rm col}(\Omega_1w_1^3, \sqrt{-1}\Omega_2 w_2, -\sqrt{-1}\Omega_2 \bar w_2)$, we obtain the parameter set
$$ \Pi_1=\left\{a\in \Pi_0: \quad |\langle k,\omega\rangle\pm m\Omega_2|\geq \frac{\gamma_0}{|k|^{n_0+1}}, m=1,2; 0<|k|\leq K_0 \right\}$$
 with the measure estimate
 \begin{equation}\label{measest0}
 {\rm Meas} \Pi_1\geq {\rm Meas}\Pi_0(1-c_3\gamma_0),
 \end{equation}
  and the estimates of the change
\begin{eqnarray}
|\| U_{1,l_100}|\|_{r_1,\Pi_1} & \lessdot &  \gamma_0^{-1}\rho_0^{-(2n_0+1)}\varepsilon_0^{3+l_1},\qquad l_1=0,1, \nonumber\\
|\| U_{1,l_100}|\|_{r_1,\Pi_1} & \lessdot & \gamma_0^{-2}\rho_0^{-(3n_0+2)}\varepsilon_0^{3+l_1},\qquad l_1=2,3, \nonumber\\
|\|U_{1,0l_2l_3}|\|_{r_1,\Pi_1} & \lessdot & \gamma_0^{-2}\rho_0^{-(3n_0+2)}\varepsilon_0^5,\qquad l_2+l_3=1, \nonumber\\
|\| U_{j,l_100}|\|_{r_1,\Pi_1} & \lessdot & \gamma_0^{-2}\rho_0^{-(3n_0+2)}\varepsilon_0^{2+l_1},\qquad l_1=0,1,\label{T0est}\\
|\| U_{j,200}|\|_{r_1,\Pi_1} & \lessdot & \gamma_0^{-3}\rho_0^{-(4n_0+3)}\varepsilon_0^4, \nonumber\\
|\| U_{j,300}|\|_{r_1,\Pi_1} & \lessdot & (\varepsilon_0+\gamma_0^{-3}\rho_0^{-(4n_0+3)}\varepsilon_0^5), \nonumber\\
|\|U_{j,0l_2l_3}|\|_{r_1,\Pi_1} & \lessdot & \gamma_0^{-2}\rho_0^{-(3n_0+2)}\varepsilon_0^4,\qquad l_2+l_3=1; j=2,3.\nonumber
\end{eqnarray}
The drifted terms from $G$ are
$$\widehat{N}(u)={\rm col}\left(\sum_{l_1=0}^3\varepsilon_0^{1+l_1} \widehat{G_{1,l_100}}(0; a,\varepsilon_0^4)u_1^{l_1}, \varepsilon_0^4\widehat{G_{2,010}}(0; a, \varepsilon_0^4)u_2, \varepsilon_0^4 \widehat{G_{3,001}}(0; a,\varepsilon_0^4)\bar u_2\right).$$
We reduce $N(u)+\widehat{N}(u)$ to a normal form as in \eqref{odev} by the translation
$$\mathscr{T}_2^0:\quad u_1=w_{10}+w_{1+},\, u_2=w_{2+},\, \bar u_2=\bar w_{2+},$$
where the $w_{10}$ is a real root of the algebraical equation
\begin{equation}\label{aeq}
\Omega_1u_1^3+\sum_{l_1=0}^3\varepsilon_0^{1+l_1} \widehat{G_{1,l_100}}(0; a,\varepsilon_0^4)u_1^{l_1}=0,
\end{equation}
and based on the assumption (H4), equivalently, the condition \eqref{case1} or \eqref{case2}, is determined in the following manner.

{\bf Case 1}  In the case where $\inf_{a\in \Pi_0}|\widehat{G_{1,000}}(0; a,0)|\geq c_1$, that is, the condition \eqref{case1} holds, the equation \eqref{aeq} has a real root $w_{10}$ satisfying
$$w_{10}=-\left(\varepsilon_0(\Omega_1^1)^{-1}\widehat{G_{1,000}}(0; a,\varepsilon_0^4 )\right)^{\frac{1}{3}}[1+O(\varepsilon_0)]$$
by Assumptions (H1) and (H3), which implies
\begin{equation}\label{3u10est}
\inf_{a\in \Pi_0}|w_{10}|\geq c_4 \varepsilon_0^{\frac{1}{3}},\qquad |\|w_{10}|\|_{\Pi_0}\leq c_5 \varepsilon_0^{\frac{1}{3}}
\end{equation}
for sufficiently small $\varepsilon_0$, where $c_4$ and $c_5$ are positive constants, $\Omega_1^1=\Omega_1+\varepsilon_0^4 \widehat{G_{1,300}}(0; a,\varepsilon_0^4)$.

{\bf Case 2} If $\widehat{G_{1,000}}(0; a,0)\equiv 0$ but $\inf_{a\in \Pi_0}|\widehat{G_{1,100}}(0; a,0)|\geq c_1$, that is, the condition \eqref{case2} holds, then let $w_{10}=0$ and put $\varepsilon_0 \widehat{G_{1,000}}(0; a,\varepsilon_0^4)(=O(\varepsilon_0^5))$ in the lower-degree terms on the next step.

Thus, the change $\mathscr{T}^0=\mathscr{T}_1^0 \circ \mathscr{T}_2^0:\, \mathscr{W}(s_1)\rightarrow \mathscr{W}(s_0-\delta_0)$  in the form
$w=w_+ +W_0^0 +W^0(\phi, w_+; a) (\phi\in \mathscr{U}(r_1), a\in \Pi_1) $,
transforms the equation \eqref{ode3}, still denoting the new variable $w_+$ by $w$, into the form
\begin{equation}\label{ode4}
 \left\{\begin{array}{rl}
\dot{w}_1 & =\varepsilon_0^2\left[N_1^1(w_1;a)+ F_1^1(\phi,w;a)+G_1^1(\phi,w; a)\right] \\
\dot{w}_2 & =\Omega_2^1(a)w_2+ F_2^1(\phi,w;a)+G_2^1(\phi,w; a)\\
\dot{\bar w}_2 & =\Omega_3^1(a)\bar w_2+ F_3^1(\phi,w;a)+G_3^1(\phi,w; a)
\end{array}
\right.
 \end{equation}
with the normal form
$$N^1={\rm col}(N_1^1(w_1;a),\Omega_2^1(a)w_2,\Omega_3^1(a)\bar w_2),$$
where $F^1_j$ and $G^1_j, j=1,2,3$ are the higher-degree and lower-degree terms, respectively,
$$N_1^1(w_1;a)=\Omega_1^1(a)w_1^3+e_2^1(a)w_1^2+e_1^1(a)w_1,\quad \Omega_1^1(a)=\Omega_1+\varepsilon_0^4\widehat{G_{1,300}}(0; a, \varepsilon_0^4), $$
 $$\Omega_2^1(a)=\sqrt{-1}\Omega_2+\varepsilon_0^4\widehat{G_{2,010}}(0; a, \varepsilon_0^4), \quad  \Omega_3^1(a)=-\sqrt{-1}\Omega_2+\varepsilon_0^4 \widehat{G_{3,001}}(0; a,\varepsilon_0^4),$$
 $$e_2^1(a)=\left\{\begin{array}{ll}
 3\Omega_1^1w_{10}^0+\varepsilon_0^3\widehat{G_{1,200}}(0; a, \varepsilon_0^4),\quad & {\rm in \, Case \, 1}\\
 \varepsilon_0^3\widehat{G_{1,200}}(0; a, \varepsilon_0^4),\quad & {\rm in \, Case \, 2,}
 \end{array}\right.$$
$$e_1^1(a)=\left\{\begin{array}{ll}
 3\Omega_1^1(w_{10}^0)^2+2\varepsilon_0^3\widehat{G_{1,200}}(0; a, \varepsilon_0^4)w_{10}^0+\varepsilon_0^2\widehat{G_{1,100}}(0; a, \varepsilon_0^4),\quad & {\rm in \, Case \, 1}\\
 \varepsilon_0^2\widehat{G_{1,100}}(0; a, \varepsilon_0^4),\quad & {\rm in \, Case \, 2.}
 \end{array}\right.$$

Now we verify that the equation \eqref{ode4} satisfies Conditions $(\nu.1)-(\nu.3)$ with $\nu=1$ in Lemma 1. By the assumption (H3) and the estimate \eqref{3u10est}, it implies that there exist a constant $c_2$ such that
$$\inf_{a\in \Pi_0}|e_1^1(a)|\geq c_2 \varepsilon_0^{\kappa},\qquad |\|e_1^1|\|_{\Pi_0}\lessdot \varepsilon_0^{\kappa},\qquad |\|e_2^1|\|_{\Pi_0}\lessdot \varepsilon_0^{\beta},$$
where, for the case 1, we take
$$\beta=\frac{1}{3}, \qquad \kappa= \frac{2}{3}, \qquad \iota=\frac{1}{2}, \qquad 0<\varepsilon_0\leq \gamma_0^{12},$$
and for the case 2
$$\beta=1, \qquad \kappa=2, \qquad \iota=\frac{7}{9}, \qquad 0<\varepsilon_0\leq \gamma_0^{9}.$$
Set
 $$\Omega_1^0=\Omega_1,  \qquad \Omega_2^0=\sqrt{-1}\Omega_2,  \qquad \Omega_3^0=-\sqrt{-1}\Omega_2,  \qquad \epsilon_1=\varepsilon_0^{\iota}.$$
Then, by the assumptions (H1) and (H3), the reality condition (H5), estimates \eqref{T0est} and \eqref{3u10est},  it is easy to see the conditions $(\nu.1)-(\nu.3)$ with $\nu=1$ in Lemma 1 hold for sufficiently small $\varepsilon_0$.

Hence, we use Lemma 1 inductively to obtain a sequence of quasi-periodic coordinate transformations $\mathscr{T}^{\nu}:\mathscr{W}_{\nu+1}\rightarrow \mathscr{W}_\nu$ in the form
$$\mathscr{T}^{\nu}:\qquad w=w_+ +W_0^\nu(a) +W^\nu(\phi, w_+; a),$$
which is analytic in $(\phi,w_+)\in \mathscr{D}(r_{\nu+1},s_{\nu+1})$ and continuously differentiable in $a\in \Pi_{\nu+1}, \nu=0,1,2,\cdots$. Set
$$\mathscr{T}_{\nu}=\mathscr{T}^0\circ\mathscr{T}^1\circ \cdots \circ\mathscr{T}^{\nu},\qquad \mathscr{D}_{\infty}= \mathscr{D}(\frac{r_0}{2}, \frac{s_0}{2}),\qquad \Pi_{\gamma_0}=\bigcap_{\nu=0}^{\infty}\Pi_\nu.$$

{\bf 1) Measure estimate}

By \eqref{measest}, \eqref{measest0} and $\Pi_{\nu+1}\subset \Pi_{\nu}(\nu=0,1,2,\cdots)$, we have
$${\rm Meas}(\Pi_0-\Pi_{\gamma_0})\leq \sum_{\nu=0}^{\infty}{\rm Meas}(\Pi_\nu-\Pi_{\nu+1})\leq c_3 {\rm Meas}\Pi_0 \sum_{\nu=0}^{\infty}\frac{\gamma_0}{(\nu+1)^2},$$
that is,
$${\rm Meas}\Pi_{\gamma_0}={\rm Meas} \Pi_0-O(\gamma_0).$$

\vskip 0.2in
{\bf 2) Convergence}

By the transformation $\mathscr{T}_{\nu}$, the equation \eqref{ode3} becomes
\begin{equation}\label{odevpl1}
 \left\{\begin{array}{rl}
\dot{w}_{1+} & =\varepsilon_0^2\left[N_1^{\nu+1}(w_{1+};a)+ F_1^{\nu+1}(\phi,w_+;a)+G_1^{\nu+1}(\phi,w_+; a)\right] \\
\dot{w}_{2+} & =\Omega_2^{\nu+1}(a)w_{2+}+ F_2^{\nu+1}(\phi,w_+;a)+G_2^{\nu+1}(\phi,w_+; a)\\
\dot{\bar w}_{2+} & =\Omega_3^{\nu+1}(a)\bar w_{2+}+ F_3^{\nu+1}(\phi,w_+;a)+G_3^{\nu+1}(\phi,w_+; a).
\end{array}
\right.
 \end{equation}
From \eqref{Jacobiest} and \eqref{w10est}-\eqref{Wjl2l3} it follows
\begin{equation}\label{convest}
|\|\mathscr{T}_{\nu}-\mathscr{T}_{\nu-1}|\|_{r_{\nu+1},s_{\nu+1},\Pi_{\nu+1}}\leq \prod_{i=0}^{\nu-1} |\|D \mathscr{T}_{i}|\|_{r_{i+1},s_{i+1},\Pi_{i+1}} |\|W_0^{\nu}+W^{\nu}|\|_{r_{\nu+1},s_{\nu+1},\Pi_{\nu+1}}\leq 2\epsilon_{\nu}^{\frac{1}{3}}.
\end{equation}
Thus, $\{\mathscr{T}_{\nu}\}$ is convergent under the norm $||\cdot||_{\frac{r_0}{2}, \frac{s_0}{2}, \Pi_{\gamma_0}}$. Let the limit be $\mathscr{T}_{\infty}$. Then $\mathscr{T}_{\infty}$ is of the form
$$\mathscr{T}^{\infty}:\qquad w=w_+ +W_0^\infty(a) +W^\infty(\phi, w_+; a),$$
the $W^\infty(\phi, w_+; a)$ is analytic in $(\phi,w_+)\in \mathscr{D}_{\infty}$ and satisfies the reality condition \eqref{rc}, $W_0^\infty(a)=(w_{10}^\infty(a),0,0)^T$ with $w_{10}^\infty(a)=\sum_{\nu=0}^{\infty}w_{10}^\nu(a)$ being real, and $W_0^\infty$ and $W^\infty$ are Lipschitz  in $a\in \Pi_{\gamma_0}$ by \eqref{convest} (the proof of "Lipschitz" is similar to the one in \cite{LL} and see the proof of Corollary 6.5 in \cite{LL}). By Lemma 1 and letting $\nu\rightarrow \infty$ in \eqref{odevpl1}, the $\mathscr{T}^{\infty}$ transforms \eqref{ode3} into
\begin{equation}\label{odeinf}
 \left\{\begin{array}{rl}
\dot{w}_{1+} & =\varepsilon_0^2\left[N_1^{\infty}(w_{1+};a)+ F_1^{\infty}(\phi,w_+;a)\right] \\
\dot{w}_{2+} & =\Omega_2^{\infty}(a)w_{2+}+ F_2^{\infty}(\phi,w_+;a)\\
\dot{\bar w}_{2+} & =\Omega_3^{\infty}(a)\bar w_{2+}+ F_3^{\infty}(\phi,w_+;a),
\end{array}
\right.
 \end{equation}
where $F_j^{\infty}$ $(j=1,2,3)$ are the higher-degree terms. Obviously, the $w_+=0$ is a solution of \eqref{odeinf}, hence $w=W_0^\infty(a) +W^\infty(\phi, 0; a)$ is a solution of \eqref{ode3}. By \eqref{rescaleT}, \eqref{complexT} and the reality condition \eqref{rc}, we obtain the quasi-periodic solution of \eqref{ode1}
$$v_1=\varepsilon^{\frac{1}{4}}(w_{10}^\infty(a)+W_1^\infty(\omega t, 0; a)), \quad v_2=\varepsilon^{\frac{1}{2}} {\rm Re}(W_2^\infty(\omega t, 0; a)), \quad
\quad v_3=\varepsilon^{\frac{1}{2}} {\rm Im}(W_3^\infty(\omega t, 0; a)),$$
which is real analytic in $\omega t$, Lipschitz in $a\in \Pi_{\gamma_0}$ and satisfies the estimate \eqref{solest} noting estimates \eqref{T0est} and \eqref{3u10est}. The proof of Theorem 1 is complete.

\renewcommand{\theequation}{\thesection.\arabic{equation}}
\section*{6. Proof of Theorem 2}
\setcounter{section}{4}\setcounter{equation}{0}

As the proof procedure of Theorem 2 is supplied in the same manner as in Theorem 1, we only describe the proof sketch. After transformed by \eqref{complexT} and \eqref{rescaleT} with $z_t\rightarrow \varepsilon^{\frac{1}{2}}z_t$, the equation \eqref{dde1} can be written in the form
\begin{equation}\label{dde2}
 \left\{\begin{array}{rl}
\dot{w}_1 & =\varepsilon_0^2\left[\Omega_1w_1^3+\varepsilon_0 F_1^0+\sum_{(l,m,n)\in \mathfrak{\Sigma}_L}\varepsilon_0^{2(|l|+|m|+|n|)+1-l_1} G_{1,lmn}(\phi)w^lz_t^m(0)z_t^n(-1)\right] \\
\dot{w}_2 & =\sqrt{-1}\Omega_2w_2+\varepsilon_0 d_4 w_1^3+\varepsilon_0^2F_2^0+ \sum_{(l,m,n)\in \mathfrak{\Sigma}_L}\varepsilon_0^{2(|l|+|m|+|n|)+2-l_1} G_{2,lmn}(\phi)w^lz_t^m(0)z_t^n(-1)\\
\dot{\bar w}_2 & =-\sqrt{-1}\Omega_2\bar w_2+\varepsilon_0 \overline{d_4}w_1^3+\varepsilon_0^2F_3^0+ \sum_{(l,m,n)\in \mathfrak{\Sigma}_L}\varepsilon_0^{2(|l|+|m|+|n|)+2-l_1} G_{3,lmn}(\phi)w^lz_t^m(0)z_t^n(-1)\\
\frac{dz_t}{dt} & =U_{Q}z_t+ X_0^Q [\varepsilon_0d_3v_1^3+\varepsilon_0^2F_4^0+\sum_{(l,m,n)\in \mathfrak{\Sigma}_L}\varepsilon_0^{2(|l|+|m|+|n|)+2-l_1} G_{4,lmn}(\phi)w^lz_t^m(0)z_t^n(-1)],
\end{array}
\right.
 \end{equation}
here we have suppressed the dependence of all functions on $a$ and $\varepsilon_0$, $F_j^0=F_j^0(\phi,w,z_t(0),z_t(-1)) (j=1,2,3)$ are the higher-degree terms, the index set of the lower-degree terms is
$$\Sigma_L=\left\{(l,m,n):\quad l_2+l_3+|m|+|n|=0\,{\rm and}\, 0\leq l_1\leq 3, \, {\rm or}\,l_2+l_3+|m|+|n|=1\,{\rm and}\, l_1=0 \right\}$$
and the coefficient $G_{1,lmn}(=G_{1,(l_1l_2l_3)mn}$ satisfies the condition \eqref{case1} or \eqref{case2}, that is
$$\inf_{a \in \Pi_0}\left|\widehat{G_{1,(000)00}}(0;a,0)\right|\geq c_1,$$
or
$$\widehat{G_{1,(000)00}}(0;a,0)=0 \quad {\rm for}\, a\in \Pi_0,\qquad {\rm and} \qquad
\inf_{a \in \Pi_0}\left|\widehat{G_{1,(100)00}}(0;a,0)\right|\geq c_1.$$
After reducing the equation \eqref{dde2} at the initial step, similar to the first part in Section 5, at the iteration step we consider the following equation
\begin{equation}\label{ddev}
 \left\{\begin{array}{rl}
\dot{w}_1 & =\varepsilon_0^2\left[N_1^\nu(w_1)+ F_1^\nu(\phi,w,z_t(0),z_t(-1))+ G_1^\nu(\phi,w,z_t(0),z_t(-1))\right] \\
\dot{w}_2 & =\Omega_2^\nu w_2+ F_2^\nu(\phi,w,z_t(0),z_t(-1))+ G_2^\nu(\phi,w,z_t(0),z_t(-1))\\
\dot{\bar w}_2 & =\Omega_3^\nu\bar w_2+ F_3^\nu(\phi,w,z_t(0),z_t(-1))+ G_3^\nu(\phi,w,z_t(0),z_t(-1))\\
\frac{dz_t}{dt} & =U_{Q}z_t+M_1^\nu(\phi,w,z_t(0),z_t(-1))+ F_5^\nu(\phi,w,z_t(0),z_t(-1))+G_5^\nu(\phi,w_1)\\
& \quad + X_0^Q [M_2^\nu(\phi,w,z_t(0),z_t(-1))+ F_4^\nu(\phi,w,z_t(0),z_t(-1))+G_4^\nu(\phi,w_1)],
\end{array}
\right.
 \end{equation}
where
$$ N_1^\nu(w_1)=\Omega^\nu_1 w_1^3+e_2^\nu w_1^2+e_1^\nu w_1,$$
$$M_i^\nu(\phi,w,z_t(0),z_t(-1))=\sum_{l_2+l_3+|m|+|n|=1} M_{i,(0l_2l_3)mn}^\nu(\phi)w_2^{l_2}\bar w_2^{l_3}z_t^m(0)z_t^n(-1), \quad i=1,2.$$
Then, we can prove Theorem 2 by using the following iteration lemma.

\vskip 0.2in
 \begin{lemma}\label{iterativelemma2}
 Assume $\epsilon_1\leq \gamma_0^6$ and that at the $\nu$-th ($\nu\geq 1$) step the equation \eqref{ddev} satisfies\\
{\bf $(\nu. 1)$ (Frequency condition)} it is the same as in Lemma 1;\\
{\bf $(\nu. 2)$ (Smallness condition)}  for $1\leq j \leq 5$ and $i=1,2$, $F_j^\nu, G_j^\nu, M_i^\nu\in \mathfrak{A}_{r_\nu,s_\nu}^{\Pi_\nu}$ ( in the present case, $\mathscr{W}(s_\nu)=\{(w,z)\in \mathbb{C}^3\times Q^{\mathbb{C}}:\, |w|\leq s_\nu, ||z||\leq s_\nu\}$, the $Q^{\mathbb{C}}$ is the complexification of $Q$), $F_j^\nu$ and $G_j^\nu$ are the higher-degree terms and the lower-degree small perturbation terms respectively, $G_4^\nu$ and $G_5^\nu$ only depend on $w_1$, $M_{1,(0l_2l_3)00}^\nu(\phi)\in Q^{\mathbb{C}}\bigcap \mathcal{C}^{1,\mathbb{C}}(l_2+l_3=1)$, $M_{1,0mn}^\nu(\phi): \mathbb{C}^2 \rightarrow Q^{\mathbb{C}}\bigcap \mathcal{C}^{1,\mathbb{C}}(|m|+|n|=1)$ (here, the $\mathcal{C}^{1,\mathbb{C}}$ is the complexification of $\mathcal{C}^1= C^1([-1,0],\mathbb{R}^q)$),  with
$$|\|F_j^\nu-F_j^{\nu-1}|\|_{r_\nu,s_\nu,\Pi_\nu}\leq \epsilon_\nu^{\frac{1}{3}},\qquad  G^\nu_{j,(l_100)00}=O_{r_\nu,\Pi_\nu}(\epsilon_\nu^{4-l_1}),\quad 0\leq l_1\leq 3, 1\leq j\leq 5 ,$$
 $$ |\|M_{i,(0l_2l_3)mn}^\nu-M_{i,(0l_2l_3)mn}^{\nu-1}|\|_{r_\nu,s_\nu,\Pi_\nu}\leq \epsilon_\nu^{\frac{1}{3}},\,  G^\nu_{j,(0l_2l_3)mn}=O_{r_\nu,\Pi_\nu}(\epsilon_\nu),\quad l_2+l_3+|m|+|n|=1, 1\leq j\leq 3, i=1,2;$$
{\bf $(\nu. 3)$ (Reality condition)}  denoting $\chi=(w_1,w_2,\bar w_2,z_t(0),z_t(-1)), \underline{\chi}=(w_1,\bar w_2,w_2,z_t(0),z_t(-1))$, the RHSEs of $\dot w_1, \dot w_2,\dot {\bar w}_2$ and $\frac{dz_t}{dt}$ in (\ref{ddev}) by $H_1^\nu,\cdots, H_4^\nu$ respectively, for $\phi \in \mathbb{T}^{n_0}$
$$ \overline{H_2^\nu}(\phi,\underline{\chi})=H_3^\nu(\phi,\chi),\quad \overline{H_3^\nu}(\phi,\underline{\chi})=H_2^\nu(\phi,\chi),\quad \overline{H_j^\nu}(\phi,\underline{\chi})=H_j^\nu(\phi,\chi),\, j=1,4.$$
Then there is a closed subset $\Pi_{\nu+1}\subset \Pi_\nu$ of the measure estimate
$${\rm Meas} \Pi_{\nu+1}\geq {\rm Meas}\Pi_\nu(1-c_3\gamma_\nu),$$
 with a constant $c_3>0$, and a quasi-periodic coordinate transformation
$\mathscr{T}^{\nu}:\mathscr{W}_{\nu+1}\rightarrow \mathscr{W}(s_\nu-\delta_\nu)\subset \mathscr{W}_{\nu}$ in the form
$$\mathscr{T}^{\nu}:\qquad w=w_+ +W_0^\nu +W^\nu(\phi, w_+,y_t(0),y_t(-1)),\quad z_t=y_t+W_4^\nu (\phi,w_{1+}),$$
where $\phi\in \mathscr{U}(r_{\nu+1}), a\in \Pi_{\nu+1}$, $w_+=(w_{1+},w_{2+},\bar w_{2+})^T$ and $y_t$ are the new variables, $W_0^\nu=(w_{10}^\nu(a),0,0)^T$,
$W^\nu=(W^\nu_1,W^\nu_2,W^\nu_3)^T, W^\nu_j \in \mathfrak{A}_{r_{\nu+1},s_{\nu+1}}^{\Pi_{\nu+1}}$ is of the same form as $G_j^\nu(j=1,2,3)$ and $W^\nu_4 \in \mathfrak{A}_{r_{\nu+1},s_{\nu+1}}^{\Pi_{\nu+1}}$ as $G_5^\nu$, satisfying Reality condition
$$w_{10}^\nu\in \mathbb{R},\quad
\overline{W_1^\nu}(\phi,\underline{\chi_+})=W_1^\nu(\phi,\chi_+),\quad \overline{W_4^\nu}(\phi,w_{1+})=W_4^\nu(\phi,w_{1+}),$$
$$\overline{W_2^\nu}(\phi,\underline{\chi_+})=W_3^\nu(\phi,\chi_+),\qquad \overline{W_3^\nu}(\phi,\underline{\chi_+})=W_2^\nu(\phi,\psi_+)$$
 for $\phi \in \mathbb{T}^{n_0}$, where $\chi_+=(w_{1+},w_{2+},\bar w_{2+},y_t(0),y_t(-1))$, and the following estimates
\begin{eqnarray*}
& & w_{10}^\nu = O_{\Pi_{\nu+1}}(\epsilon_{\nu}^{\frac{4}{3}}), \quad
W_{1,(l_100)00}^\nu = O_{r_{\nu+1},\Pi_{\nu+1}}(\varepsilon_0^2\epsilon_{\nu}^{4-l_1}\gamma_{\nu}^{-2}\rho_{\nu}^{-(3n_0+2)}), 0\leq l_1\leq 3,\\
& & W_{1,(0l_2l_3)00}^\nu = O_{r_{\nu+1},\Pi_{\nu+1}}(\varepsilon_0^2\epsilon_{\nu}\gamma_{\nu}^{-2}\rho_{\nu}^{-(3n_0+2)}),\quad l_2+l_3=1,\\
& & W_{j,(l_100)00}^\nu = O_{r_{\nu+1},\Pi_{\nu+1}}(\epsilon_{\nu}^{4-l_1}\gamma_{\nu}^{-2}\rho_{\nu}^{-(3n_0+2)}),\quad 0\leq l_1\leq 3; j=2,3,4,\\
& & W_{j,(0l_2l_3)00}^\nu = O_{r_{\nu+1},\Pi_{\nu+1}}(\epsilon_{\nu}\gamma_{\nu}^{-2}\rho_{\nu}^{-(3n_0+2)}),\quad l_2+l_3=1, j=2,3,\\
& &  W_{1,0mn}^\nu = O_{r_{\nu+1},\Pi_{\nu+1}}(\varepsilon_0^2\epsilon_{\nu}),\quad , W_{j,0mn}^\nu = O_{r_{\nu+1},\Pi_{\nu+1}}(\epsilon_{\nu}),j=2,3;\quad |m|+|n|=1,
\end{eqnarray*}
such that the equation (\ref{ddev}) is transformed into the same form satisfying Conditions $(\nu.1)-(\nu.3)$ by replacing $\nu$ with $\nu+1$ and $(w_1,w_2,\bar w_2,z_t)$ with $(w_{1+},w_{2+},\bar w_{2+},y_t)$, respectively.
\end{lemma}

{\bf Outline of proof}  The proof is similar to that of Lemma1, only adding arguments of the hyperbolic part. We reduce the lower-degree terms only depending on $w_1$ in the equation of the hyperbolic direction. Drop the index $\nu$ and rewrite the equation \eqref{ddev} as
\begin{equation}\label{dde3}
\mathfrak{J}_{\varepsilon_0}^{-1}\left(\begin{array}{c}\dot w\\ \frac{dz_t}{dt}\end{array}\right) =\left(\begin{array}{c}N(w)\\ U_Q z_t \end{array}\right)
+ \left(\begin{array}{c}G+F\\ M_1+G_5+F_5+X_0^Q(M_2+G_4+F_4)\end{array}\right),
 \end{equation}
where, $\mathfrak{J}_{\varepsilon_0}={\rm diag}(J_{\varepsilon_0}, {\rm Id})$.

Take the first transformation as the following form
$$ \mathscr{T}_1:\qquad w=u+U(\phi,u,y_t(0),y_t(-1)),\quad z_t=y_t+U_4(\phi,u_1)$$
with $u=(u_1,u_2,\bar u_2)^T$. We introduce some notations. For a given function $z(\theta)$ from $[-1,0]$ to $\mathbb{C}^q$, set
$$U_{j,01}(\phi)z=U_{j,0m0}(\phi)z^m(0)+U_{j,00n}(\phi)z^n(-1),\, G_{j,01}(\phi)z=G_{j,0m0}(\phi)z^m(0)+G_{j,00n}(\phi)z^n(-1),\, 1\leq j\leq 3,$$
$$M_{i,01}(\phi)z=M_{i,0m0}(\phi)z^m(0)+M_{i,00n}(\phi)z^n(-1),\quad i=1,2; \qquad |m|+|n|=1.$$

Substituting the transformation $\mathscr{T}_1$ into (\ref{dde3}), we obtain the homological equations
\begin{eqnarray}
 & &\omega \partial_\phi U_{1,(l_100)00}+(l_1-1)\varepsilon_0^2e_1U_{1,(l_100)00} + (l_1-3)\varepsilon_0^2e_2U_{1,((l_1-1)00)00}\nonumber\\
 & & \quad +(l_1-5)\varepsilon_0^2\Omega_1U_{1,((l_1-2)00)00}= \varepsilon_0^2(\Gamma_K G_{1,(l_100)00}(\phi)-\widehat{G_{1,(l_100)00}}(0)),\, 0\leq l_1\leq 3,\label{4U1l1}\\
& & \omega \partial_\phi U_{1,(0l_2l_3)00}-\varepsilon_0^2e_1U_{1,(0l_2l_3)00} +  l_2\Omega_2U_{1,(0l_2l_3)00} +l_3\Omega_3 U_{1,(0l_2l_3)00} \nonumber\\
&  &  =\Gamma_K [\varepsilon_0^2 G_{1,(0l_2l_3)00}(\phi)-U_{1,01}(M_{1,(0l_2l_3)00}(\phi)+X_0^QM_{2,(0l_2l_3)00}(\phi))],\quad l_2+l_3=1,\label{4U1l2l3}\\
& & \omega \partial_\phi U_{1,01}-\varepsilon_0^2e_1U_{1,01} + U_{1,01}[U_Q+M_{1,01}(\phi)+X_0^QM_{2,01}(\phi)]
 =\varepsilon_0^2 G_{1,01}(\phi),\label{4U1mn}\\
 & & \omega \partial_\phi U_{j,(l_100)00}-\Omega_j U_{j,(l_100)00}+l_1\varepsilon_0^2e_1U_{j,(l_100)00} + (l_1-1)\varepsilon_0^2e_2U_{j,((l_1-1)00)00} \nonumber\\
& & \quad +(l_1-2)\varepsilon_0^2\Omega_1U_{j,((l_1-2)00)00}= \Gamma_K G_{j,(l_100)00}(\phi),\quad 0\leq l_1\leq 3; j=2,3,\label{4U23l1}\\
& & \omega \partial_\phi U_{2,(010)00}=\Gamma_K \widetilde{G}_{2,(010)00}(\phi)-\widehat{\widetilde{G}_{2,(010)00}}(0),\label{4U2l2}\\
 & & \omega \partial_\phi U_{2,(001)00}+(\Omega_3-\Omega_2)U_{2,(001)00}= \Gamma_K\widetilde{ G}_{2,(001)00}(\phi),\label{4U2l3}\\
 & & \omega \partial_\phi U_{3,(010)00}-(\Omega_3-\Omega_2)U_{3,(010)00}= \Gamma_K \widetilde{G}_{3,(010)00}(\phi), \label{4U3l2}\\
& & \omega \partial_\phi U_{3,(001)00}=\Gamma_K \widetilde{G}_{3,(001)00}(\phi)-\widehat{\widetilde{G}_{3,(001)00}}(0),\label{4U3l3}\\
& & \omega \partial_\phi U_{j,01}-\Omega_j U_{j,01} + U_{j,01}[U_Q+M_{1,01}(\phi)+X_0^QM_{2,01}(\phi)]
 =G_{j,01}(\phi), \quad j=2,3, \label{4U23mn}\\
& & \omega \partial_\phi U_{4,(l_100)00}-[U_Q+M_{1,01}(\phi)+X_0^QM_{2,01}(\phi)]U_{4,(l_100)00} +l_1\varepsilon_0^2e_1U_{4,(l_100)00} \nonumber\\
& & \quad + (l_1-1)\varepsilon_0^2e_2U_{4,((l_1-1)00)00} +(l_1-2)\varepsilon_0^2\Omega_1U_{4,((l_1-2)00)00} \nonumber\\
 & & = G_{5,(l_100)00}(\phi)+X_0^Q G_{4,(l_100)00}(\phi)+  [M_{1,(010)00}(\phi)+X_0^QM_{2,(010)00}(\phi)]U_{2,(l_100)00} \nonumber\\
& & \quad + [M_{1,(001)00}(\phi)+X_0^QM_{2,(001)00}(\phi)]U_{3,(l_100)00},\quad 0\leq l_1\leq 3, \label{4U4l1}
\end{eqnarray}
where, $U_{j,(l_100)00}=0$ if $l_1<0$ for $1\leq j \leq 4$,
$$\widetilde{G}_{j,(0l_2l_3)00}(\phi)=G_{j,(0l_2l_3)00}(\phi)-U_{j,01}(\phi)[M_{1,(0l_2l_3)00}(\phi)+X_0^QM_{2,(0l_2l_3)00}(\phi)],\quad l_2+l_3=1, j=2,3.$$
We first solve the equations \eqref{4U1l1}, \eqref{4U1mn}, \eqref{4U23l1} and \eqref{4U23mn}, then \eqref{4U1l2l3}, \eqref{4U2l2}-\eqref{4U3l3} and \eqref{4U4l1}. These equations \eqref{4U1mn}, \eqref{4U23mn} and \eqref{4U4l1} do not involve small divisors and are solved by the same method as for (4.12) and (4.14) in \cite{LY12}, the other equations by the same manner as in Lemma 1. The remainder of the proof of the lemma is similar to that of Lemma 1 and the detail is omitted.

\vskip 0.2in
\appendix
\renewcommand{\theequation}{\thesection.\arabic{equation}}
\section*{ Appendix }
\setcounter{section}{1}\setcounter{equation}{0}

The following lemma shows that the reality condition is preserved by the transformation of a Jacobian matrix with the reality condition or its inverse matrix. For the sake of simplicity, we introduce some notations. Let $\mathfrak{X} $ denote the set of all vector-functions $F(\phi,w)={\rm col}(F_1(\phi,w), F_2(\phi,w), F_3(\phi,w))$, being analytic in $(\phi,w)\in \mathscr{U}(r)\times \mathscr{W}(s)$ and satisfying the reality condition
\begin{equation}\label{A1}
F_1(\phi,w)=\overline{F_1}(\phi,\underline{w}),\quad  F_2(\phi,w)=\overline{F_3}(\phi,\underline{w}),\quad  F_3(\phi,w)=\overline{F_2}(\phi,\underline{w})
\end{equation}
for $\phi \in \mathbb{T}^{n_0}$, where $w=(w_1,w_2,\bar w_2)^T$ and $\underline{w}=(w_1,\bar w_2,w_2)^T$, and $\frac{DF}{Dw}(\phi,w)$ represent the Jacobian matrix $\frac{D(F_1,F_2,F_3)}{D(w_1,w_2,\bar w_2)}$ of $F$ with respect to $w$.  Define $\underline{F}\coloneqq {\rm col}(F_1(\phi,w), F_3(\phi,w), F_2(\phi,w))$. The reality condition means $\overline{F}(\phi,\underline{w})=\underline{F}(\phi,w)$ for $\phi \in \mathbb{T}^{n_0}$.

\vskip 0.2in
 \begin{lemma}\label{reality}
 Assume $F(\phi,w), H(\phi,w)\in \mathfrak{X}$. Then $\left(\frac{DH}{Dw} F\right)(\phi,w)\in \mathfrak{X}$. Moreover, if the Jacobian matrix $\frac{DH}{Dw}(\phi,w)$ is invertible for all $(\phi,w)\in \mathscr{U}(r)\times \mathscr{W}(s)$, then $\left(\frac{DH}{Dw}\right)^{-1} F(\phi,w)\in \mathfrak{X}$.
 \end{lemma}

 {\bf Proof} Set $B(\phi,w)=\left(\frac{DH}{Dw} F\right)(\phi,w)$. If $F, H\in \mathfrak{X}$, the reality condition \eqref{A1} implies
 \begin{eqnarray}
 \overline{B}(\phi,\underline{w}) & = & \overline{\left(\frac{DH}{Dw}\right)}(\phi,\underline{w}) \overline{F}(\phi,\underline{w}) \nonumber\\
  & = & \left(\frac{D\underline{H}}{D\underline{w}}\right)(\phi,w) \overline{F}(\phi,\underline{w})=\left(\frac{D\underline{H}}{D\underline{w}}\right)(\phi,w) \underline{F}(\phi,w) \label{A2}\\
  & = & \underline{\left(\frac{DH}{Dw} F\right)}(\phi,w)=\underline{B}(\phi,w),\nonumber
  \end{eqnarray}
  that is, $B(\phi,w)\in \mathfrak{X}$. Moreover, if $H\in \mathfrak{X}$ and the $\frac{DH}{Dw}$ is invertible, from $\left(\frac{DH}{Dw} F\right)(\phi,w)\in \mathfrak{X}$ and the second equality in \eqref{A2}, it follows the second part of the lemma.

\vskip 0.2in
 \begin{lemma}\label{root}
 Let $f(x;a)=b_3(a)x^3+b_2(a)x^2+b_1(a)x+b_0(a)$ be a real third-order polynomial in $x$. Suppose the coefficients $b_j(a)(0\leq j \leq 3)$ are continuously differentiable in $a\in \Pi$ ($\Pi \subset \mathbb{R}$ is a bounded closed set), and there exist positive constants $c_j(6\leq j \leq 8)$ such that
 $$|\|b_0|\|_{\Pi}\lessdot \epsilon^4, \qquad |\|b_1|\|_{\Pi}\lessdot \varepsilon_0^\kappa, \qquad
 \inf_{a\in \Pi}|b_1(a)|\geq c_6 \varepsilon_0^\kappa,$$
 $$|\|b_2|\|_{\Pi}\lessdot \varepsilon_0^\beta, \qquad |\|b_3|\|_{\Pi}\leq c_7, \qquad   \inf_{a\in \Pi}|b_3(a)|\geq c_8,$$
 where
 $$\epsilon \leq \varepsilon_0^\iota, \qquad \beta <\frac{4}{3}\iota, \qquad \kappa <\beta +\frac{4}{3}\iota.$$
 Then the equation $f(x;a)=0$ has a real root $x_0(a)$ which is continuously differentiable with respect to the parameter $a\in \Pi$ and satisfies the estimate
 $|\|x_0|\|_{\Pi}\leq c_9 \epsilon^{\frac{4}{3}}$, here $c_9$ is a positive constant.
 \end{lemma}

 {\bf Proof} By Theorem 2.1 in \cite{Xu} it is obvious that the equation $f(x;a)=0$ has a real root $x_0(a)$ such that
 \begin{equation}\label{x0est}
 \sup_{a\in \Pi}|x_0(a)| \lessdot \epsilon^{\frac{4}{3}}.
 \end{equation}
 Noting that the assumption of the lemma implies
 $$|\partial_x f(x_0;a)|=|3b_3(a)x_0^2+2b_2(a)x_0+b_1(a)|\geq \frac{c_6}{2}\varepsilon_0^\kappa>0$$
 for sufficiently small $\varepsilon_0$, and $f(x;a)$ is continuously differentiable in $a\in \Pi$, we  deduce that $x_0(a)$ is also continuously differentiable in $a\in \Pi$ and satisfies
 $$\sup_{a\in \Pi}\left|\frac{d}{da}x_0(a)\right| \lessdot \epsilon^{\frac{4}{3}}$$
 based on differentiating the equation $f(x_0;a)\equiv 0$ with respect to $a$, and using \eqref{x0est} and the assumption. The proof is complete.


\vskip 0.3in

\begin{thebibliography}{99}
\bibitem{Ar63} V. I. Arnol'd, Small denominators and problems of stability of motion in classical mechanics and celestial mechanics, Uspekhi Mat. Nauk, 18(1963) 91-192.
\bibitem{BBM11} D. Bambusi, M. Benti, E. Magistrelli, Degenerate KAM theory for partial differential equations, J. Differential Equations, 250(2011) 3379-3397.
\bibitem{BG01} D. Bambusi, G. Gaeta, Invariant tori for non-conservative perturbations of integrable systems, NoDEA Nonlinear Differential Equations Appl., 8(2001) 99-116.
\bibitem{BB87} B. J. L. Braaksma and H. W. Broer, On a quasi-periodic Hopf bifurcation, Ann. Inst. H. Poincar\'{e}, Anal. nonlin., 4(1987) 115-168.
 \bibitem{BBH90} B. J. L. Braaksma, H. W. Broer and G. B. Huitema, Toward a quasi-periodic bifurcation theory, Mem. Amer. Math. Soc., 83(421)(1990) 83-167.
\bibitem{BDL14} J. Bramburger, B. Dionne, V. G. LeBlanc, Zero-Hopf bifurcation in the van der Pol oscillator with delayed position and velocity feedback, Nonlinear Dyn., 78(2014) 2959-2973.
\bibitem{BHS96} H. W. Broer, G. B. Huitema, M. B. Sevryuk, Quasi-periodic Motions in Families of Dynamical Systems: Order amidst Chaos, Lecture Notes in Math., Vol.1645, Springer, Berlin, 1996.
\bibitem{Bru92} A. D. Bruno, On conditions for nondegeneracy in Kolmogorov's theorem, Soviet Math. Dokl., 45(1992) 221-225.
\bibitem{Ch96} C.-Q. Cheng, Birkhoff-Kolmogorov-Arnold -Moser tori in convex Hamiltonian systems, Comm. Math. Phys., 177(1996) 529-559.
\bibitem{Ch99} C.-Q. Cheng, Lower-dimensional invariant tori in the regions of instability for nearly integrable Hamiltonian systems, Comm. Math. Phys., 203(1999) 385-419.
\bibitem{CS94} C.-Q. Cheng, Y. Sun, Existence of KAM tori in degenerate Hamiltonian systems, J. Differential Equations, 114(1994) 288-335.
\bibitem{CW99} C.-Q. Cheng, S. Wang, The surviving of lower dimensional tori from a resonant torus of Hamiltonian systems, J. Differential Equations, 155(1999) 311-326.
\bibitem{Fr67} M. Friedman, Quasi-periodic solutions of nonlinear ordinary differential equations with small damping, Bull. Amer. Math. Soc., 73(1967) 460-464.
\bibitem{Gen07} G. Gentile, Degenerate lower-dimensional tori under the Bryuno condition, Ergod. Th. \& Dynam. Sys., 27(2007) 427-457.
\bibitem{GG05} G. Gentile, G. Gallavotti, Degenerate elliptic resonances, Comm. Math. Phys., 257(2005) 319-362.
\bibitem{HLY10} Y. Han, Y. Li, Y. Yi, Invariant  tori in Hamiltonian systems with high order proper degeneracy, Ann. Henri Poincar\'{e}, 10(2010) 1419-1436.
\bibitem{HLY} Y. Han, Y. Li, Y. Yi, Degenerate lower-dimensional tori in Hamiltonian systems, J. Differential Equations, 227(2006) 670-691.
\bibitem{JW08} W. Jiang, J. Wei, Bifurcation analysis in van der Pol's oscillator with delayed feedback, J. Comp. Appl. Math., 213(2008) 604-615.
\bibitem{JS} A. Jorba, C. Simo, On quasi-periodic perturbations of elliptic equilibrium points, SIAM J. Math. Anal., 27(1996) 1704-1737.
\bibitem{Li16} X. Li, On the persistence of quasi-periodic invariant tori for double Hopf bifurcation of vector fields, J. Differential Equations, 260(2016) 7320-7357.
\bibitem{LL} X. Li,  R. de la Llave, Construction of quasi-periodic solutions of delay differential equations via KAM technique, J. Differential Equations, 247(2009) 822-865.
\bibitem{LY12} X. Li, X. Yuan, Quasi-periodic solutions for perturbed autonomous delay differential equations, J. Differential Equations, 252(2012) 3752-3796.
\bibitem{LY03} Y. Li, Y. Yi, A quasi-periodic Poincare's Theorem, Math. Ann., 326(2003) 649-690.
\bibitem{Mos67} J. Moser, Convergent series expansions for quasi-periodic motions, Math. Ann., 169(1967) 136-176.
\bibitem{Ru01} H. R\"{u}ssmann, Invariant tori in non-degenerate nearly integrable Hamiltonian systems, Regul. Chaotic Dyn., 6(2001) 119-204.
\bibitem{Sev07} M. B. Sevryuk, Invariant tori in quasi-periodic non-autonomous dynamical systems via Herman's method, Discrete Contin. Dyn. Syst., 18(2007) 569-595.
\bibitem{SM71} C. L. Siegel, J. Moser, Lectures on Celestial Mechanics, Springer, Berlin, 1971.
\bibitem{WJ10} H. Wang, W. Jiang, Hopf-pitchfork bifurcation in van der Pol's oscillator with nonlinear delay feedback, J. Math. Anal. Appl., 368(2010) 9-18.
\bibitem{Xu} J. Xu, On small perturbation of two-dimensional quasi-periodic systems with hyperbolic-type degenerate equilibrium point, J. Differential Equations, 250(2011) 551-571.
\bibitem{XC03} J. Xu, K. W. Chung, Effects of time delayed position feedback on a van der Pol-Duffing oscillator, Physica D, 180(2003) 17-39.
\bibitem{You} J. You, A KAM theorem for hyperbolic-type degenerate lower dimensional tori in Hamiltonian systems, Comm. Math. Phys., 192(1998) 145-168.
\bibitem{Yuan02} X. Yuan, Construction of quasi-periodic breathers via KAM technique, Commum. Math. Phys., 226(2002) 61-100.
\bibitem{ZG13} L. Zhang, S. Guo, Hopf bifurcation in delayed van der Pol oscillators, Nonlinear Dyn., 71(2013) 555-568.



\end{thebibliography}
\end{document}